\documentclass[aos,MSNbibl,seceqn,dvips]{arximspdf}
\makeatletter
   \@ifpackageloaded{graphicx}{}{\usepackage{graphicx}}
\makeatother
\usepackage{dcolumn}

\doi{10.1214/15-AOS1347}% Updated by VTEXPTS2LaTeX.exe, 08.07.2015
\volume{43}
\issue{6}
\pubyear{2015}
\firstpage{2451}
\lastpage{2483}
\docsubty{FLA}

\makeatletter
\newcolumntype{d}[1]{D{.}{.}{#1}}
\newcommand{\rrvert}{\vert}
\newcommand{\rrVert}{\Vert}
\newcommand{\llvert}{\vert}
\newcommand{\llVert}{\Vert}
\renewcommand{\mid}{|}
\newcommand{\xrightarrow}[1]{\stackrel{#1}{\rightarrow}}
\newcommand{\eqref}[1]{(\ref{#1})}
\newcommand{\mathds}{\mathbb}

\newtheorem{teo}{Theorem}[section]
\newtheorem{cor}[teo]{Corollary}
\newtheorem{prop}[teo]{Proposition}
\newproclaim{rem}[teo]{Remark}
\newproclaim{ex}[teo]{Example}
\newproclaim{ass}[teo]{Assumption}
\newproclaim{alg}[teo]{Algorithm}

\newcommand{\argmax}{\operatorname{argmax}}

\newcommand{\xiv}{\bolds{\xi}}
\newcommand{\muv}{\bolds{\mu}}
\newcommand{\ind}{\mathbf{1}}
\newcommand{\R}{\mathds{R}}
\newcommand{\Z}{\mathds{Z}}
\newcommand{\N}{\mathds{N}}

\newcommand{\cov}{\operatorname{Cov}} %\mathbb{C}\textnormal{o\hspace*{0.02cm}v}}
\newcommand{\E}{\mathbb{E}}
\newcommand{\W}{\mathcal{W}}
\newcommand{\OO}{\mathcal{O}}
\newcommand{\Gaussian}{\mathcal{N}}

\renewcommand{\ss}{s}
\newcommand{\ad}{\mathfrak{a}}
\newcommand{\bd}{\mathfrak{b}}
\newcommand{\ud}{\mathfrak{u}}
\newcommand{\dd}{\mathfrak{d}}
\newcommand{\ld}{\mathfrak{l}}
\newcommand{\pd}{\mathfrak{p}}
\newcommand{\qd}{\mathfrak{q}}
\newcommand{\td}{\mathfrak{t}}
\newcommand{\B}{\mathcal{B}}
\newcommand{\bs}[1]{\bolds{#1}}
\newcommand{\WW}{\bs{\mathcal{W}}}
\newcommand{\BB}{\bs{\mathcal{B}}}
\renewcommand{\ss}{s}
\newcommand{\Tone}{\textup{(T1)}}
\newcommand{\Ttwo}{\textup{(T2)}}
\newcommand{\Sone}{\textup{(S1)}}
\newcommand{\Stwo}{\textup{(S2)}}
\newcommand{\Sthree}{\textup{(S3)}}
\newcommand{\Bone}{\textup{(B1)}}
\newcommand{\Btwo}{\textup{(B2)}}
\newcommand{\Mboot}{\textup{(M)}}
\newcommand{\SRboot}{\textup{(SR)}}
\newcommand{\SNRboot}{\textup{(SNR)}}
\makeatother

\begin{document}
\begin{frontmatter}

\title{Uniform change point tests in high dimension}
\runtitle{Uniform change point tests}

\begin{aug}
% Corresponding author: Moritz Jirak - m0ritz@yahoo.com% Updated by
%VTEXPTS2LaTeX.exe, 08.07.2015 15:19
\author[A]{\fnms{Moritz}~\snm{Jirak}\corref{}\thanksref{T2}\ead[label=e1]{jirak@math.hu-berlin.de}}
\runauthor{M. Jirak}
\affiliation{Humboldt Universit\"{a}t zu Berlin}
%\dedicated{}
\address[A]{Institut f\"{u}r Mathematik\\
Unter den Linden 6\\
D-10099 Berlin\\
Germany\\
\printead{e1}}
\end{aug}
\thankstext{T2}{Supported by the Deutsche Forschungsgemeinschaft via
FOR 1735
\textit{Structural Inference in Statistics}: \textit{Adaptation and
Efficiency}.}

% HISTORY:
%
\received{\smonth{8} \syear{2014}}% Updated by VTEXPTS2LaTeX.exe,
%08.07.2015 15:19
%
\revised{\smonth{4} \syear{2015}}% Updated by VTEXPTS2LaTeX.exe,
%08.07.2015 15:19

% ABSTRACT
%
\begin{abstract}
Consider $d$ dependent change point tests, each based on a
CUSUM-statistic. We provide an asymptotic theory that
allows us to deal with the maximum over all test statistics as both the
sample size $n$ and $d$ tend to infinity. We achieve this either by a
consistent bootstrap or an appropriate limit distribution. This allows
for the construction of simultaneous confidence bands for dependent
change point tests, and explicitly allows us to determine the location of
the change both in time and coordinates in high-dimensional time
series. If the underlying data has sample size greater or equal $n$ for
each test, our conditions explicitly allow for the large $d$ small $n$
situation, that is, where $n/d \to0$. The setup for the high-dimensional time series is based on a general weak dependence concept.
The conditions are very flexible and include many popular multivariate
linear and nonlinear models from the literature, such as ARMA, GARCH
and related models. The construction of the tests is completely
nonparametric, difficulties associated with parametric model selection,
model fitting and parameter estimation are avoided. Among other things,
the limit distribution for $\max_{1 \leq h \leq d} \sup_{0 \leq t
\leq1}\llvert  \mathcal{W}_{t,h} - t \mathcal{W}_{1,h}\rrvert  $ is
established, where $ \{\mathcal{W}_{t,h} \}_{1 \leq h \leq
d}$ denotes a sequence of dependent Brownian motions. As an
application, we analyze all S\&P 500 companies over a period of one year.
\end{abstract}

% KEYWORDS
\end{frontmatter}

%s1 #&#
\section{Introduction}\label{sec1}

Modeling high-dimensional time series is a necessity in many different
fields, ranging from meteorological and agricultural problems to
biology, genetics, financial engineering and risk management. %, see
%for instance~\cite{Diebold2003,hallinforni2000,lee1992,%nelsonsiegel1987,worsley2002}.
Particularly within the financial regulation framework, banks and
insurance undertakings are required to assess and incorporate hundreds
of different factors and risks. Regarding financial time series, it is
well known that large panels of asset returns routinely display break
points and other nonstationarities (cf.~\cite{fan2011}). In this
context, structural stability is a very important issue, since even
changes in few parameters can lead to misspecified risk measures and
wrong conclusions (cf.~\cite{esma,spokoiny2009}). The issue of
structural stability also arises in many other fields, such as
climatology, genetics and medicine. Hence, given a $d$-dimensional time
series ${\mathbf X}_k =  (X_{k,1},\ldots,X_{k,d} )^{\top}$,
there is a high interest in procedures that consistently partition the
coordinates of $ \{{\mathbf X}_k \}_{k \in\Z}$ into the two sets
\begin{eqnarray}
\nonumber
\mathcal{S}_d &=& \bigl\{1 \leq h \leq d:
\{X_{k,h}\}_{k \in\Z}\mbox { is stable} \bigr\},
\nonumber\\[-8pt]\\[-8pt]\nonumber
\mathcal{S}_d^c &=& \bigl\{1 \leq h \leq d:
\{X_{k,h}\}_{k \in\Z}\mbox{ is unstable} \bigr\},
\end{eqnarray}
such that we have the relation
\begin{eqnarray*}
\{1,\ldots,d \} &=& \mathcal{S}_d \uplus\mathcal{S}_d^c.
\end{eqnarray*}
The sample included in $\mathcal{S}_d$ may then be used for further
inference, while the part contained in $\mathcal{S}_d^c$ requires
subsequent treatment. Often, one is additionally interested in the
actual time of change in each coordinate, and tests based on cumulative
sums are efficient in this context. Let us denote such tests with
$B_{n,h}^{\widehat{\sigma}}$, $1 \leq h \leq d$ for further
reference; see \eqref{defnBnh} below for a precise definition.

In the univariate case, tests for structural stability in time series
are widely available (cf.~\cite
{auehorvathsurvey2012,banerjeeetal2005,berkesgombayhi,csoghovbook1993,changepointhorvatcsogobook1997,perron2005}
and the many references there), the multivariate setup, and especially
the high-dimensional case are less often considered. Apart from
functional data approaches (cf.~\cite
{berkes2009JRRS,fremdt2013,hoermann2010}), the literature in the
latter case is rather sparse compared to the univariate theory. Let us
briefly mention some recent contributions in this area. In~\cite
{horvath2012}, the stability of panel data is considered. Using a
threshold-aggregation approach,~\cite{ChoFryzlewicz2012Preprint}
study the detection of global changes (see also~\cite
{groen2010,ombao2005}), whereas in~\cite{Jirak2012chcov}, the
topic of possible gain or loss in power in higher dimension is
discussed. Changes in the covariance structure in a multivariate setup
are addressed in~\cite{aueetal2009}, and an interesting connection
between Dos-attacks and change point detection is explored in~\cite
{leduc2009} (see also~\cite{tartakovskychange2004,xie2013} and therein for changes in multi-channel systems). However,
to the best of my knowledge, a (thorough) treatment regarding the
consistent estimation of $\mathcal{S}_d^c$, particularly in a time
series framework, is lacking in the literature so far. Compared to the
univariate case, handling the multivariate situation is much more
complicated since breaks may or may not be present at different times
in different coordinates $h$. Since it is usually unknown which
coordinates $h$ have anomalies and which ones have not, determining
$\mathcal{S}_d^c$ (resp., $\mathcal{S}_d$) is particularly hard if
the dimension $d$ is large. The vast majority of high-dimensional
change point procedures use aggregation or PCA based techniques, and
are therefore inappropriate for determining $\mathcal{S}_d^c$. In this
context, a natural way to measure possible deviations is to employ the
statistic $T_d^{\widehat{\sigma}} = \max_{1\leq h \leq d}
B_{n,h}^{\widehat{\sigma}}$, with coordinate-wise CUSUM-statistics
%
%e1.1 #&#
\begin{equation}
\label{defnBnh} \qquad B_{n,h}^{\widehat{\sigma}} = \bigl(\widehat{
\sigma}_h^2 n\bigr)^{-1/2}\max
_{1 \leq k \leq n}\Biggl\llvert \sum_{j = 1}^{k}
X_{j,h} - \frac{k}{n} \sum_{j = 1}^n
X_{j,h}\Biggr\rrvert, \qquad h = 1,\ldots,d.
\end{equation}
Here, $\widehat{\sigma}_h^2$ is an appropriate estimator for the
long-run variance, which will be more fully explained below. Control of
$T_d^{\widehat{\sigma}}$ readily allows us to make inference for every
single coordinate $h$. In this paper, we provide theoretic tools that
allows one to handle $T_d^{\widehat{\sigma}}$. If the random
variables $X_{k_1,h_1}$ and $X_{k_2,h_2}$ become less dependent if
either quantity $\llvert  k_1-k_2\rrvert  $ or $\llvert  h_1-h_2\rrvert  $ becomes large, we will show that
%
%e1.2 #&#
\begin{equation}
\label{eqintroextremecoonvcusum} \max_{1 \leq h \leq d} e_d
\bigl(B_{n,h}^{\widehat{\sigma}} - f_d \bigr) \xrightarrow{w}
\mathcal{V} \qquad\mbox{as both $n,d \to \infty$,}
\end{equation}
for appropriate normalizing sequences $e_d, f_d$, where $\mathcal{V}$
is an extreme value distribution of Gumbel type. A general explicit
connection between $n$ and $d = d_n$ is given such that \eqref
{eqintroextremecoonvcusum} is valid for $n,d \to\infty$, allowing
for $n/d \to0$, but also for the converse where $d/n \to C \geq0$. On
the other hand, we show that the time series $\{X_{k,h}\}_{k \in\Z}$
may have properties such that a pivotal limit theorem like in \eqref
{eqintroextremecoonvcusum} cannot exist. For this case, we provide
bootstrap approximations, which of course work in both cases.

Studying the joint limit as $d,n \to\infty$ is a much more realistic
setup than considering its sequential analogue (i.e., $\lim_{d \to
\infty} \lim_{n \to\infty} \cdot$, cf. Remark 2.1 in~\cite
{aueetal2009} or~\cite{groen2010}), but is also considerably
harder from a mathematical point of view. In order to allow for a high
flexibility in \eqref{eqintroextremecoonvcusum}, we use a
generalization of known weak dependence concepts from the univariate
(multivariate) case to the high-dimensional setup, which allows for
dependencies in time and space. This leads to fairly general, yet
easily verifiable conditions that are valid for a large number of
popular time series from the literature, including multivariate ARMA
and GARCH models. Even though we only consider breaks in the mean
vector, it is clear that our results are also applicable for assessing
the stability of the variance or second-order structure (possibly
cross-wise) up to a certain extent.

An outline of the paper can be given as follows. In Section~\ref
{secmain}, we introduce and discuss our assumptions and the main
results. The aspect of concise estimation of $\mathcal{S}_d^c$ and the
actual time of change within $\mathcal{S}_d^c$ is discussed in
Section~\ref{sectime}. Bootstrap procedures and their consistency are
explored in Section~\ref{secpermalt}. Section~\ref{secexamples}
contains a number of popular time series examples that are included in
our framework. Section~\ref{secempresults} deals with practical
aspects and investigates the finite sample behavior. As a real data
application, we simultaneously analyze all S\&P 500 companies over the
time horizon of one year. %This data set contains some of the a fore
%mentioned problems and characteristicae.
Detailed proofs are given in the supplementary material \cite{suppA}.

%s2 #&#
\section{Methodology and main results}\label{secmain}

Throughout this paper, we use $\lesssim$, $\gtrsim$, ($\thicksim$)
to denote (two-sided) inequalities involving a multiplicative constant.
$C$ denotes an arbitrary, absolute constant that may vary from line to
line. Let $\llVert  \cdot\rrVert  _p$ denote the $\mathds{L}^p$-norm
$\E[|\cdot
|^p]^{1/p}$ for $p \geq1$, and given a set $\mathcal{S}$, we write
$\llvert  \mathcal{S}\rrvert  $ to symbolize its cardinality. We write $\stackrel
{d}{=}$ for equality in distribution. In the sequel, we often deal with
arrays $ (c_h )_{1 \leq h \leq d}$, where $d \to\infty$ and
$c_h$ may depend on $d$. We then use the abbreviations
%
%e2.1 #&#
\begin{equation}
\label{defninfsupstar} \inf_h^* c_h = \liminf
_{d \to\infty} \min_{1 \leq h \leq d} c_h, \qquad
\sup_h^* c_h = \limsup_{d \to\infty}
\max_{1 \leq h \leq d} c_h.
\end{equation}

Let\vspace*{2pt} $ \{{\mathbf X}_k \}_{k \in\Z}$ with ${\mathbf X}_k =
(X_{k,1},\ldots,X_{k,d} )^{\top}$ be a sequence of
$d$-dimensional random vectors where $\E [{\mathbf X}_k ] = \muv
_k =  (\mu_{k,1},\ldots,\mu_{k,d} )^{\top}$. The aim of
this paper is to provide a simultaneous test for structural stability
in $\muv_k$, based on the observations ${\mathbf X}_1,\ldots,{\mathbf X}_n$.
To do so, we consider the coordinate-wise null-hypothesis
%
%e2.2 #&#
\begin{equation}
\mathcal{H}_{0,h}: \mu_{1,h} = \cdots= \mu_{n,h},
\qquad h = 1,\ldots,d,
\end{equation}
which indicates structural stability in the mean over time. Under this
notion of stability, we get that $\mathcal{S}_d =  \{1 \leq h
\leq d:   \mathcal{H}_{0,h}\mbox{ is true} \}$, that is,
$\mathcal{S}_d$ denotes the set of all coordinates where $\mathcal
{H}_{0,h}$ holds. We say that $\mathcal{H}_0$ is true, if $\mathcal
{S}_d =  \{1,\ldots,d \}$. As alternative hypothesis, we
specify the scenario which allows for at most one change in each
coordinate of $\muv_k$. More precisely, we assume that there exists a
(usually unknown) time lag $k_h^*$, such that
%
%e2.3 #&#
%e2.4 #&#
\begin{eqnarray}
\mathcal{H}_{A,h}: \mu_{1,h} = \cdots= \mu_{k_h^*,h}
\neq\mu _{k_h^* + 1,h} = \cdots= \mu_{n,h}
\nonumber\\[-8pt]\\[-8pt]
\eqntext{\mbox{for some }h =
1,\ldots,d.}
\end{eqnarray}
We say that the alternative $\mathcal{H}_A$ holds, if at least one
$\mathcal{H}_{A,h}$ is true, and the null hypothesis $\mathcal
{H}_{0,h}$ hold in all the remaining unaffected coordinates. This means
there is at least one break in one coordinate $h$. Generalizations to
multiple change point detection are possible, but will not be addressed
here. We make the following convention. The Type I error refers to a
 ``false alarm,'' that is, a break detection where there is none, and
the Type II error is attributed to an unreported break. In this spirit,
we then obtain $\mathcal{S}_d^c =  \{1 \leq h \leq d:   \mathcal
{H}_{A,h}\mbox{ is true} \}$, that is, the set which consists of
all coordinates where a change has occurred. In order to identify
$\mathcal{S}_d^c$, we propose to use the coordinate-wise\vspace*{1pt} CUSUM
statistic $ B_{n,h}^{\widehat{\sigma}}$ defined in \eqref
{defnBnh}. We denote the whole vector of such statistics with ${\mathbf
B}_{n,d}^{\widehat{\sigma}} =  (B_{n,1}^{\widehat{\sigma
}},\ldots,B_{n,d}^{\widehat{\sigma}} )^{\top}$. Let
%
%e2.5 #&#
\begin{equation}
\label{eqdefnBH} \B_h = \sup_{0\leq t \leq1}\llvert
\W_{t,h} - t \W_{1,h}\rrvert, \qquad h = 1,\ldots,d,
\end{equation}
where $\WW_{t,d} =  (\W_{t,1},\ldots,\W_{t,d} )^{\top}$
is a $d$-dimensional Brownian motion, with correlations $\rho_{i,j} =
\E [\W_{t,i} \W_{t,j} ]$. If the dimension $d$ is fixed
and $n \to\infty$, it is known that under quite general conditions we
have weak convergence, that is,
\[
{\mathbf B}_{n,d}^{\widehat{\sigma}} \xrightarrow{w} \BB_d,
\]
where $\BB_d =  (\B_{1},\ldots,\B_{d} )^{\top}$, with
associated correlation matrix ${\bolds\Sigma}_d =\break  (\rho
_{i,j} )_{1 \leq i,j \leq d}$. Given some mild regularity
conditions for ${\bolds\Sigma}_d$, we will show in Theorem~A.2 in~\cite{suppA} that
%
%e2.6 #&#
\begin{equation}
\label{eqextremecoonvBB} \max_{1 \leq h \leq d} e_d (
\mathcal{B}_h - f_d ) \xrightarrow{w} \mathcal{V},
\qquad\mbox{as $d \to\infty$,}
\end{equation}
for appropriate sequences $e_d, f_d$, where $\mathcal{V}$ is an
extreme value distribution of Gumbel type. Result \eqref
{eqextremecoonvBB} is one of the key ingredients in our proof for
\eqref{eqintroextremecoonvcusum}, and may be of independent
interest. Limit theorems involving the maximum of partial sums have
played a fundamental role in statistic and probability theory for a
long time (cf.~\cite{shorak1986}). Particularly the seminal
contribution in~\cite{darling1956} has stimulated much research in
this area; see, for instance,~\cite
{csoghovbook1993,changepointhorvatcsogobook1997} for an account
on further developments and applications, and~\cite
{csorgo2003erdoes} for some sharp results and a brief historic
review. Related research can also be found in~\cite{leadbetter1974},
see also the references therein.

Based on an asymptotic result like \eqref{eqextremecoonvBB},
simultaneous confidence regions can readily be constructed, we refer to
\eqref{eqconfregion} for more details. However, non-Gaussianity is
often the rule rather than the exception. It is therefore of
considerable interest to formulate our results in a more general
manner. In the univariate case, a highly accepted model in the
literature is to assume the structure $X_k = g (\varepsilon
_k,\varepsilon_{k-1},\ldots )$ for a process $ \{X_k \}
_{k \in\Z}$, where $ \{\varepsilon_k \}_{k \in\Z}$ is a
sequence of i.i.d. random variables in some space $\mathbb{S}$ of
possible infinite dimension. Let $ \{\varepsilon_k' \}_{k \in
\Z}$ be an independent copy of $ \{\varepsilon_k \}_{k \in\Z
}$. Then many well-known weak-dependence measures and concepts are
based on quantifying the difference (for $p \geq1$)
%
%e2.7 #&#
\begin{equation}
\label{defnakuni} \qquad a_k(p) = \bigl\llVert g (\varepsilon_k,
\varepsilon_{k-1},\ldots,\varepsilon _0,\varepsilon_{-1},
\ldots ) - g \bigl(\varepsilon_k,\varepsilon _{k-1},\ldots,
\varepsilon_0',\varepsilon_{-1},\ldots \bigr)\bigr
\rrVert _p.
\end{equation}
For example, the dependence concept in~\cite{wu2005} is based on
$a_k(p)$. In related cases, the whole past is replaced with copies; see~\cite{berkes2007,poetscher1999} and~\cite{aueetal2009} for a
multivariate version. We will see in Section~\ref{secexamples} that
many well-known univariate and multivariate time series such as ARMA
and GARCH-models are within this framework. As is outlined, for example,
in~\cite{aueetal2009}, such conditions have several advantages over
certain mixing competitors. For instance, mixing conditions are
sometimes hard to verify and may require additional smoothness
assumptions (cf.~\cite{andrews1984}). A more profound discussion is
given in~\cite{wu2005}. Another advantage is that these dependence
measures have a natural spatial extension which includes the univariate
(multivariate) case as a special example; see, for instance,~\cite
{ElMachkouri2013,chenwuaos2013}. More precisely, for
$ \{{\mathbf X}_k \}_{k \in\Z}$ with ${\mathbf X}_k =  \{
X_{k,h} \}_{h \in\N}$ we have the structure condition
%
%e2.8 #&#
\begin{equation}
\label{eqstruccondiX} X_{k,h} = g_h (\varepsilon_{k},
\varepsilon_{k-1},\ldots ), \qquad k \in\Z, h \in\N,
\end{equation}
where $g_h$ are measurable functions. The coordinate processes
$X_{k,h}$ can be viewed as projections from $\mathbb{S}$ to $\R$. In
analogy to \eqref{defnakuni}, for $p \geq1$ we put [recall $\sup_h^*$ in~\eqref{defninfsupstar}]
%
%e2.9 #&#
\begin{eqnarray}\label{defnakhigh}
a_k(p) &=& \sup_h^*\bigl\llVert
g_h (\varepsilon_k,\varepsilon_{k-1},\ldots,
\varepsilon_0,\varepsilon_{-1},\ldots )
\nonumber\\[-8pt]\\[-8pt]\nonumber
&&{}  - g_h \bigl(
\varepsilon _k,\varepsilon_{k-1},\ldots,\varepsilon_0',
\varepsilon_{-1},\ldots \bigr)\bigr\rrVert _p.
\end{eqnarray}
Note that $a_k(p)$ is a \textit{temporal} dependence measure, that is,
it only measures dependence in time, and essentially does not impose
any \textit{spatial} dependence restrictions. As extreme possibly
examples just consider the cases where $X_{k,h} = X_{k,h+1}$ are
identical or where $ \{X_{k,h} \}_{1 \leq h \leq d}$ is an
independent sequence for each $k \in\Z$. In fact, this setup is very
general and contains a huge variety of popular linear and nonlinear
time series models, see Section~\ref{secexamples} for more details.

Allowing for weak dependence in (multivariate) time series inevitably
results in dealing with the long run covariances
$\gamma_{i,j}$, which we formally introduce as
%
%e2.10 #&#
\begin{equation}
\label{eqvardef} \gamma_{i,j} = \lim_{n \to\infty}
n^{-1}\E \Biggl[\sum_{k =
1}^{n}
\sum_{l = 1}^n (X_{k,i} -
\mu_{k,i}) (X_{l,j} - \mu _{l,j}) \Biggr].
\end{equation}
We shall see (cf.~\cite{suppA}) that Assumption~\ref{assmain} below
implies that the above limit exists and $\gamma_{i,j}$ are thus well
defined. Moreover, in case of $\sigma_h^2 \stackrel{\mathrm{def}}{=} \gamma
_{h,h}$ we have the usual representation $\sigma_h^2 = \sum_{k \in\Z
} \phi_{k,h}$, where $\phi_{k,h} = \cov[X_{0,h},X_{k,h}]$. If
$\sigma_{i}, \sigma_j > 0$, we also have $\rho_{i,j} = \gamma
_{i,j}/ (\sigma_i \sigma_j )$. Our main temporal assumption
is now as follows.

%as2.1 #&#
\begin{ass}[(Temporal assumptions)]\label{assmain}
Given representation \eqref{eqstruccondiX} for $\{{\mathbf X}_k\}_{k
\in\Z}$, assume that for $p > 4$ and absolute constant $\sigma^- > 0$:
\begin{longlist}
\item[\Tone]\label{A1} $a_k(p) \lesssim k^{-\ad}$, with $\ad> 5/2$,
\item[\Ttwo]\label{A2} $\inf_{h}^*\sigma_h \geq\sigma^- > 0$.
%\item[\Tthree]\label{A3}
\end{longlist}
\end{ass}

Let us briefly discuss these assumptions. \hyperref[A1]{\Tone} is a
global, polynomial decay assumption on the temporal dependence. In the
univariate case, $\ad> 1$ is possible and essentially optimal. Here,
we require the slightly stronger condition $\ad> 5/2$, which enables
us to operate in a high-dimensional context. Assumption \hyperref
[A2]{\Ttwo} is a nondegeneracy assumption that we require since we
often normalize with $\sigma_h$ in the sequel; see, however, Remark
\ref{remremoveindexset}. Note that we require Assumption~\ref
{assmain} \textit{throughout} the remainder of this paper.

Since $\sigma_h^2$ is usually unknown, we need to estimate it. The
literature (cf.~\cite{timeseriesbrockwell}) provides many potential
candidates to estimate $\sigma_{h}^2$. A popular estimator is
Bartlett's estimator, or more general, estimators of the form
%
%e2.11 #&#
\begin{equation}
\label{eqvarestimator} \widehat{\sigma}_{h}^2 = \sum
_{| k|  \leq b_n} \omega(k/b_n) \widehat {
\phi}_{k,h}, \qquad\mbox{$b_n \to\infty$,}
\end{equation}
with weight function $\omega(x)$, where
\begin{eqnarray*}
\widehat{\phi}_{h,j} &=& (n-j)^{-1}\sum
_{k = j + 1}^n (X_{k,h} - \overline{X}_{h})
(X_{k-j,h} - \overline{X}_{h}),
\end{eqnarray*}
and $\overline{X}_h = n^{-1}\sum_{k = 1}^n X_{k,h}$. Setting $\omega
(x) = 1$, we obtain the plain estimate (cf.~\cite{wu2007anals}). For
the sake of simplicity, we just consider the plain estimate for our
theoretical analysis, but the results remain equally valid for other
weight functions. Conditions on the possible size of the bandwidth $b_n
\thicksim n^{\bd}$, $0 <\bd< 1$ in terms of $\bd$ are given below in
Assumption~\ref{assrelations}.

In order to establish a limit theory, we also require some spatial
dependence conditions. A very general way that leads to easily
verifiable conditions is in terms of decay assumptions for the
underlying covariance structure. This is a common approach in the
literature; see, for instance,~\cite{cai2011,leadbetterbook1983,xiaowuspa2013}. In our context, it is
thus natural to impose conditions on the correlations $\rho_{i,j}$. As
stated in the \hyperref[sec1]{Introduction}, we consider the situation where we allow
that both $n$ and $d$ jointly tend to infinity. For a more formal
description, we model the dimension $d$ as $d \thicksim n^{\dd}$, $\dd
> 0$ throughout the remainder of this paper. The necessary connection
between $\bd, \dd$ and the underlying moments $p$ is now additionally
collected in our \textit{spatial} assumptions.

%
%as2.2 #&#
\begin{ass}[(Spatial assumptions)]\label{assrelations}
Assume that $\dd, \bd$, $ (\rho_{i,j} )_{1 \leq i,j \leq
d}$, $p > 4$ satisfy the conditions below, uniformly in $d$ for
absolute constants $\rho^+, C_{\rho},\delta> 0$:
\begin{longlist}
\item[\Sone]\label{S1} $0<\dd< \min \{p/2 - 2, (1 - \bd)p/2
-1 - \bd\}$,
\item[\Stwo]\label{S2} $\sup_{i,j: \llvert  i-j\rrvert   \geq1} \rho_{i,j} \leq
\rho^+ < 1$,
\item[\Sthree]\label{S3} $\llvert  \rho_{i,j}\rrvert   \leq C_{\rho} \log
(\llvert  i-j\rrvert  +2)^{-2-\delta}$.
\end{longlist}
\end{ass}

%
%re2.3 #&#
\begin{rem}\label{remexpgrowthconditions}
Assumptions \hyperref[S2]{\Stwo}, \hyperref[S3]{\Sthree} are only
needed for establishing the asymptotic distribution in Theorem~\ref
{teoextreme} below. Also note that the polynomial growth rate of the
dimension $d = d_n$ and the polynomial decay rate of $a_k(p)$ are
intimately connected. In this spirit, one may show that analogue
results as presented below are valid for an exponentially growing
dimension $d$, by imposing exponential decay rates on $a_k(p)$. Such
results would require in addition that $\sup_h^*\E [e^{s_0
X_{k,h}} ] <\infty$ for some $s_0 > 0$.
\end{rem}

%
%re2.4 #&#
\begin{rem}\label{remdimension}
For ease of exposition, we distinctly asked for $\dd> 0$ in \hyperref
[S1]{\Sone} to ensure that $d \to\infty$ as $n \to\infty$, which
results in a minimal polynomial growth rate. However, we point out that
we actually only require that $d \lesssim n^{\dd}$ and $d \to\infty$
as $n \to\infty$, which is slightly more general.
\end{rem}

Assumption~\ref{assrelations} only imposes mild conditions,
essentially allowing for any polynomial growth rate of the dimension $d
\thicksim n^{\dd}$ given sufficiently many moments. Note that high
moment assumptions are common in such a context, we refer to~\cite
{cai2011,jiang2004max,jirak2012,xiaowuspa2013}, where sometimes
up to 30 moments and more are required. Also note that we only need a
logarithmic decay for the correlations $\rho_{i,j}$ that is close to
the best-known results in the literature in a different context (cf.~\cite{leadbetter1974}). We are now ready to state our first main
result, which establishes the asymptotic limit distribution.

%
%th2.5 #&#
\begin{teo}\label{teoextreme}
Assume that $\mathcal{H}_0$ and Assumptions~\ref{assmain} and~\ref
{assrelations} hold. Then
\begin{eqnarray*}
\lim_{n \to\infty} P \Bigl( \max_{1 \leq h \leq d}
B_{n,h}^{\widehat{\sigma}} \leq u_d \bigl(e^{-x}
\bigr) \Bigr) = \exp\bigl(-e^{-x}\bigr),
\end{eqnarray*}
where $u_d (e^{-x} ) = x/e_d + f_d $, with $e_d = 2\sqrt
{2\log(2 d)}$, $f_d = e_d/4$.
\end{teo}

%
%re2.6 #&#
\begin{rem}\label{remremoveindexset}
Conditions \hyperref[S2]{\Stwo}, \hyperref[S3]{\Sthree} are needed
to exclude any pathologies. However, as is known in the literature (cf.~\cite{csorgo2003erdoes,csorgo2003self}), condition $\sigma_h >
0$ can be removed to some extent by more detailed arguments due to the
self-normalization in $B_{n,h}^{\widehat{\sigma}}$. Moreover, let
${\mathcal I}_d \subset \{1,\ldots,d \}$ be any sequence of
subsets with cardinality $\llvert  \mathcal{I}_d\rrvert  /d \to0$. It is then shown
in~\cite{suppA} that it actually suffices to have \hyperref[S2]{\Stwo
}, \hyperref[S3]{\Sthree} only for $i,j \in \{1,\ldots,d
\} \setminus{\mathcal I}_d$.
\end{rem}

In Section~\ref{secempresults}, we give a brief account on the
implications and relevance of the necessary assumptions for real data
sets. A problem that can appear in practical applications is the rate
of convergence to extreme value distributions, see Section~\ref
{secempresults} for details. One way out are bootstrap methods. We
first present a (comparatively) fast and easy to implement method for a
parametric bootstrap. To this end, let $ \{Z_{k,h} \}_{k \in
\Z,h \in\N}$ be a standard Gaussian IID sequence. Denote with
\begin{eqnarray}
\nonumber
\label{eqmaxdefn} B_{n,h}^{Z} &=& \frac{1}{\sqrt{n}}\max
_{1 \leq k \leq n}\Biggl\llvert \sum_{j = 1}^{k}
Z_{j,h} - \frac{k}{n}\sum_{j = 1}^n
Z_{j,h} \Biggr\rrvert\quad\mbox{and}
\nonumber\\[-8pt]\\[-8pt]\nonumber
T_{d}^{Z} &=& \max_{1 \leq h \leq d}
B_{n,h}^{Z}\quad \mbox{and recall} \quad
T_{d}^{\widehat{\sigma}} = \max_{1 \leq h
\leq d}
B_{n,h}^{\widehat{\sigma}}.
\end{eqnarray}
Next, we introduce the exact quantile $u_d^{Z}(z)$, defined as
\begin{eqnarray*}
P \bigl(B_{n,h}^{Z} \leq u_d^{Z}(z)
\bigr) = 1 - \frac{z}{d}.
\end{eqnarray*}
It then comes as no surprise that we have the following result.

%pr2.7 #&#
\begin{prop}\label{propsimpleboot}
Grant the assumptions of Theorem~\ref{teoextreme}. Then
\begin{eqnarray*}
&&\sup_{x \in\R}\bigl\llvert P \bigl(T_{d}^{Z}
\leq u_d^{Z} \bigl(e^{-x} \bigr) \bigr) - P
\bigl(T_{d}^{\widehat{\sigma}} \leq u_d^{Z}
\bigl(e^{-x} \bigr) \bigr)\bigr\rrvert = \mbox{\scriptsize $
\mathcal{O}$} (1 )\qquad\mbox{as }n \to\infty.
\end{eqnarray*}
\end{prop}

We thus obtain a very simple bootstrap method, which just requires the
generation of i.i.d. Gaussian random variables. Note that unlike to
$u_d(z)$, the quantiles $u_d^{Z}(z)$ are highly nonlinear, which seems
to make them less attractive. In practice though, it turns out that
$u_d^{Z}(z)$ often yields much better results than $u_d(z)$, also for
dependent time series. For more details and empirical results, see
Section~\ref{secempresults}. Based on Theorem~\ref{teoextreme} and
Proposition~\ref{propsimpleboot}, we can construct asymptotic honest
$1-\alpha$ confidence regions $\widehat{\mathcal{S}}_{d}(\alpha)$
and $\widehat{\mathcal{S}}_{d}^Z(\alpha)$ via
\begin{eqnarray}
\label{eqconfregion}
\nonumber
\widehat{\mathcal{S}}_{d}(\alpha) &=& \bigl\{1
\leq h \leq d: B_{n,h}^{\widehat{\sigma}} \leq x_{\alpha}/e_d
+ f_d \bigr\}, \qquad x_{\alpha} = -\log\bigl(-\log(1 -
\alpha)\bigr),\hspace*{-30pt}
\nonumber\\[-8pt]\\[-8pt]\nonumber
\widehat{\mathcal{S}}_{d}^Z(\alpha) &=& \bigl\{1 \leq h
\leq d: B_{n,h}^{\widehat{\sigma}} \leq u_d^Z(z_{\alpha})
\bigr\}, \qquad z_{\alpha} = d \bigl(1 - (1 - \alpha)^{1/d}
\bigr).
\end{eqnarray}

Let us now turn to the important question when we have less spatial
structure. As is demonstrated in Example~\ref{exfactor}, a pivotal
limit result like in Theorem~\ref{teoextreme} cannot exist if we drop
condition $\hyperref[S3]{\Sthree}$. Fortunately, things do not go
totally wrong. Our next result essentially implies that the  ``rate''
(resp., normalization) $u_d (\cdot )$ in Theorem~\ref
{teoextreme} acts as an upper bound, even under considerably less
assumptions. This is important for statistical applications, since we
remain in control of the Type I error.

%
%th2.8 #&#
\begin{teo}\label{teoupperbound}
Assume that $\mathcal{H}_0$, Assumption~\ref{assmain} and $\hyperref
[S1]{\Sone}$ hold. Then
\begin{eqnarray*}
\lim_{n \to\infty} P \Bigl( \max_{1 \leq h \leq d}
B_{n,h}^{\widehat{\sigma}} \leq u_d \bigl(e^{-x}
\bigr) \Bigr) \geq 1 - e^{-x},
\end{eqnarray*}
where $u_d (\cdot )$ is as in Theorem~\ref{teoextreme}.
\end{teo}

A useful implication of Theorem~\ref{teoupperbound} is that under
less assumptions, $\widehat{\mathcal{S}}_{d}(\alpha)$ and $\widehat
{\mathcal{S}}_{d}^Z(\alpha)$ can be modified to also supply
(asymptotic) honest confidence sets. Indeed, using the power series of
$\log(1 - \alpha)$, we obtain that
\begin{eqnarray*}
1 - \alpha\geq1 - \exp(-z_{\alpha}) = 1 - \alpha- \frac{\alpha
^2}{2} + \OO
\bigl(\alpha^3 \bigr).
\end{eqnarray*}
In particular, if we select $a = 1 - \exp(-\alpha)$, then we can
construct the confidence sets $\widehat{\mathcal{S}}_{d}(a),\widehat
{\mathcal{S}}_{d}^Z(a)$, which according to Theorem~\ref
{teoupperbound} at least have nominal level $\alpha$, since we have
[with $x_a = -\log(-\log(1 - a))$]
\begin{eqnarray*}
1 - \exp(-x_a) = 1 + \log(1 - a) = 1 - \alpha.
\end{eqnarray*}
Hence, the resulting confidence regions might be too large, but never
too small, which implies that the Type I error of the null-hypothesis
$\mathcal{H}_0$ remains controlled. Note, however, that such a
modification is more conservative, and thus results in a loss in power.
Some further properties of $\widehat{\mathcal{S}}_{d}(\alpha),
\widehat{\mathcal{S}}_{d}^Z(\alpha)$ and their behavior under the
alternative hypothesis $\mathcal{H}_A$ are the topic of Section~\ref
{sectime}. Another option to construct confidence regions if Theorem
\ref{teoextreme} fails to hold is bootstrapping. In the context of
dependent data, blockwise bootstrap procedures are a possible way out.
This topic is more fully explored in Section~\ref{secpermalt}.

%s3 #&#
\section{Estimating the location of change and general consistency of
long run variance estimation} \label{sectime}

We first make the following convention. We say that an estimator
$\widetilde{\mathcal{S}}_d$ is consistent, if
\begin{eqnarray*}
\lim_{n} P \bigl(\llvert \widetilde{\mathcal{S}}_d
\bigtriangleup \mathcal{S}_d\rrvert = 0 \bigr) = 1,
\end{eqnarray*}
where $\widetilde{\mathcal{S}}_d \bigtriangleup\mathcal{S}_d$
stands for the symmetric difference of the sets $\widetilde{\mathcal
{S}}_d$ and ${\mathcal{S}}_d$. Note that trivially any consistent
estimator $\widetilde{\mathcal{S}}_d$ gives a consistent estimator
$\widetilde{\mathcal{S}}_d^c$ for the complement set. For further
analysis, we assume that under $\mathcal{H}_A$ the times of change
depend on $n$. This is a common assumption in the literature, and one
way to guarantee this is by demanding that $k_h^* = \lfloor\tau_h n
\rfloor$ for $\tau_h \in(0,1)$. If there is no change in coordinate
$h$, we set $\tau_h = 1$. Another important quantity is the actual
minimal size of the change, which we denote with
%
%e3.1 #&#
\begin{eqnarray}
\label{defnDeltamu} \Delta\mu= \min_{h \in\mathcal{S}_d^c}\Delta\mu_h
\qquad\mbox {where }\Delta\mu_h = \llvert \mu_{k_h^*,h} -
\mu_{k_h^* + 1,h}\rrvert.
\end{eqnarray}
We assume throughout that $\Delta\mu$ is a monotone decreasing
sequence, and express the direct connection to $\mathcal{H}_A$ through
the notation $\mathcal{H}_A^{(\Delta\mu)}$. This means that under
$\mathcal{H}_A^{(\Delta\mu)}$, the minimal size of change is $\Delta
\mu$. Suppose now that $h \in\mathcal{S}_d^c$. Then elementary
calculations show that
%
%e3.2 #&#
\begin{eqnarray}
\label{eqelemtarylowerboundchange} B_{n,h}^{\widehat{\sigma}} &\geq& \widehat{
\sigma}_h^{-1}\sqrt {n}\Delta\mu_h
\tau_h (1 - \tau_h) \bigl(1 - \mbox {\scriptsize$
\mathcal{O}$}(1) \bigr)- \overline{B}_{n,h}^{\widehat
{\sigma}},
\end{eqnarray}
where
%
%e3.3 #&#
\begin{eqnarray}
\overline{B}_{n,h}^{\widehat{\sigma}} &=& \frac{1}{\widehat{\sigma
}_h\sqrt{n}}\max
_{1 \leq k \leq n}\Biggl\llvert \sum_{j = 1}^{k}
U_{j,h} - \frac{k}{n}\sum_{j = 1}^n
U_{j,h} \Biggr\rrvert,
\nonumber\\[-8pt]\\[-8pt]\nonumber
U_{j,h} &=& X_{j,h} -
\E [X_{j,h} ].
\end{eqnarray}
Due to Theorem~\ref{teoextreme}, we can control $\max_{1 \leq h \leq
d}\overline{B}_{n,h}^{\widehat{\sigma}}$ as long as $\widehat
{\sigma}_h$ behaves  ``reasonably''\vspace*{1pt} under $\mathcal{H}_A^{(\Delta\mu
)}$. If this is indeed the case, then we can expect from \eqref
{eqelemtarylowerboundchange} that $B_{n,h}^{\widehat{\sigma}}$
becomes large, and thus detect a change in coordinate $h$ using the
confidence sets $\widehat{\mathcal{S}}_{d}(\alpha)$ or $\widehat
{\mathcal{S}}_{d}^Z(\alpha)$ in \eqref{eqconfregion}.
Unfortunately though, $\widehat{\sigma}_h$ may not at all behave
reasonably and can cause the problem of  ``none monotone power'' (i.e.,
the power can decrease when the alternative gets farther away from the
null); see, for instance,~\cite{crainiceanuvogelsang2007}. One way
to overcome this problem is to use self-normalization, as proposed in~\cite{shaozhang2010}. Here, we propose a different method that will
lead to no loss in power. To this end, we first discuss the estimation
of the possible time of change $\tau_h$ for each affected coordinate.
For this, we propose the following estimates. Pick any fixed $0 < \td
\leq1/2$, preferably small, and consider
%
%e3.4 #&#
\begin{eqnarray}
\label{defntauh} \qquad \widehat{\tau}_h(\td) &=& \mathop{\argmax}_{t \in(\td,1-\td)}
\Biggl(n^{-1/2}\Biggl\llvert \sum_{j = 1}^{\lceil n t \rceil}
X_{j,h} - \frac
{\lceil n t \rceil}{n} \sum_{j = 1}^n
X_{j,h}\Biggr\rrvert \Biggr),
\qquad h = 1,\ldots,d.
\end{eqnarray}
Note that $\widehat{\sigma}_h$ is not included in the above
definition. Of course, one would like to select $\td$ such that
%
%e3.5 #&#
\begin{eqnarray}
\label{assboundtimechange} \td< \inf_{h \in\mathcal{S}_d^c}\tau_h \leq\sup
_{h \in\mathcal
{S}_d^c} \tau_h < 1 - \td.
\end{eqnarray}
In the sequel, we put $\widehat{\tau}_h = \widehat{\tau}_h(\td)$
to lighten the notation if the dependence on $\td$ is of no relevance.
Bounds for uniform deviation probabilities for $\widehat{\tau}_h(\td
)$ follow from the next result.

%th3.1 #&#
\begin{teo}\label{teotime}
Assume that $\mathcal{H}_A^{(\Delta_{\mu})}$ is valid and Assumption
\ref{assmain} holds. If we have in addition \eqref
{assboundtimechange} for some $0 < \td\leq1/2$, then
\begin{eqnarray*}
P \Bigl(\max_{h \in{\mathcal{S}}_{d}^c}\bigl\llvert \widehat{\tau
}_h(\td) - \tau_h \bigr\rrvert \geq x \Bigr) \lesssim
\bigl\llvert \mathcal {S}_d^c\bigr\rrvert (x n \log n
)^{-p/2 + 2},
\end{eqnarray*}
where we require that $x \geq C_a \frac{\log n}{(\Delta\mu)^2 n}$,
$C_a > 0$ sufficiently large.
\end{teo}

%
%re3.2 #&#
\begin{rem}
The above constant $C_a$ depends on $\td$ and the sequence $
(a_j(p) )_{j \in\N}$, and thus also implicitly on the long run
variances $ (\sigma_h )_{1 \leq h \leq d}$. Assumption $\td
> 0$ and \hyperref[A1]{\Tone} ensure that $C_a < \infty$, uniformly
in $d$.
\end{rem}

Since $\llvert  \mathcal{S}_d^c\rrvert   \leq d$, we get the following uniform
consistency result.

%co3.3 #&#
\begin{cor}\label{cortime}
Grant the assumptions of Theorem~\ref{teotime}. If in addition
\begin{eqnarray*}
\limsup_{n \to\infty} \frac{\log n}{(\Delta\mu)^2 n} = 0 \quad \mbox{and} \quad
\dd< p/2 -2,
\end{eqnarray*}
then
\begin{eqnarray*}
\max_{h \in{\mathcal{S}}_{d}^c}\bigl\llvert \widehat{\tau}_h(\td) -
\tau_h \bigr\rrvert = \mbox{\scriptsize$\mathcal{O}$}_P
(1 ).
\end{eqnarray*}
\end{cor}

Armed with Corollary~\ref{cortime}, we construct the following simple
estimators $\widehat{\sigma}_h^*$ for~$\sigma_h$:
\begin{itemize}
\item Choose a constant $0 < B_{\tau} < 1$, and use the estimator
$\widehat{\tau}_h(\td)$ in every coordinate to split the sample into
$\mathcal{T}_h^- =  \{k \leq B_{\tau} \widehat{\tau}_h(\td) n
 \}$ and $\mathcal{T}_h^+ =  \{n - B_{\tau} (1
-\widehat{\tau}_h(\td) )n < k \leq n \}$.
\item Use the usual estimator $\widehat{\sigma}_h$ to construct the
estimators $\widehat{\sigma}_h^-$ and $\widehat{\sigma}_h^+$ based
on the samples $\mathcal{T}_h^-$ and $\mathcal{T}_h^+$. The final
estimator $\widehat{\sigma}_h^*$ is then obtained by the convex combination
\begin{eqnarray*}
\bigl(\widehat{\sigma}_h^*\bigr)^2 = \widehat{
\tau}_h(\td) \bigl(\widehat{\sigma }_h^-
\bigr)^2 + \bigl(1 - \widehat{\tau}_h(\td) \bigr) \bigl(
\widehat{\sigma }_h^+\bigr)^2, \qquad1 \leq h \leq d.
\end{eqnarray*}
\end{itemize}

%
%re3.4 #&#
\begin{rem}
As was pointed out by a reviewer, possible alternatives are
\begin{eqnarray*}
\widehat{\sigma}_h^{\min} &=& \min \bigl\{\widehat{\sigma
}_h^-,\widehat{\sigma}_h^+ \bigr\}, \qquad\widehat{
\sigma}_h^{\max
} = \max \bigl\{\widehat{\sigma}_h^-,
\widehat{\sigma}_h^+ \bigr\}\quad\mbox{and}
\\
\bigl(\widehat{\sigma}_h^{\mathrm{mean}}
\bigr)^2 &=&   \bigl(\bigl(\widehat{\sigma}_h^-
\bigr)^2 + \bigl(\widehat{\sigma}_h^+\bigr)^2
\bigr)/2.
\end{eqnarray*}
Another alternative is $\widehat{\sigma}_h^{\diamond} = \widehat
{\sigma}_h^-$ if $\llvert  \mathcal{T}^-\rrvert   \geq\llvert  \mathcal{T}^+\rrvert  $ and
$\widehat{\sigma}_h^{\diamond} = \widehat{\sigma}_h^+$ otherwise.
Note that we have the relation
\begin{eqnarray*}
\widehat{\sigma}_h^{\min} \leq\widehat{
\sigma}_h^*, \widehat {\sigma}_h^{\diamond},
\widehat{\sigma}_h^{\mathrm{mean}}, \widehat{\sigma}_h^-,
\widehat{\sigma}_h^+ \leq\widehat{\sigma }_h^{\max}.
\end{eqnarray*}
As follows from the result below, all estimators are consistent, both
under $\mathcal{H}_0$ and $\mathcal{H}_A^{(\Delta_{\mu})}$, and
thus yield the correct limit distribution. From a practical perspective
though, $\widehat{\sigma}_h^{\min}$ leads to a more liberal test
with more power, whereas $\widehat{\sigma}_h^{\max}$ leads to a more
conservative test.
\end{rem}

The following result establishes the desired properties of the above
variance estimators.

%pr3.5 #&#
\begin{prop}\label{propoptimvarest}
Theorems~\ref{teoextreme} and~\ref{teoupperbound} remain valid if
we replace $\widehat{\sigma}_h$ with either $\widehat{\sigma
}_h^{\min}, \widehat{\sigma}_h^*, \widehat{\sigma}_h^{\diamond},
\widehat{\sigma}_h^{\mathrm{mean}}$, $\widehat{\sigma}_h^-$,
$\widehat{\sigma}_h^+$ or $\widehat{\sigma}_h^{\max}$. Moreover,\vspace*{1pt}
if Corollary~\ref{cortime} holds, these estimates are consistent
under $\mathcal{H}_A^{(\Delta_{\mu})}$.
\end{prop}

Now that we have settled the problem of the long run variance
estimation, we may return to our original problem of determining
$\mathcal{S}_d$, which is now an easy task. In fact, combining
Proposition~\ref{propoptimvarest} with the lower bound in \eqref
{eqelemtarylowerboundchange} and Corollary~\ref{cortime}, we
immediately get the following result.

%
%pr3.6 #&#
\begin{prop}\label{propconsistentestimationofSd}
Let $\alpha= \alpha_n \to0$ such that $\alpha_n \gtrsim n^{-1}$.
Assume in addition that \eqref{assboundtimechange} holds for some
$0 < \td_0 \leq1/2$, and that
%
%e3.6 #&#
\begin{eqnarray}
\label{eqthmtestoptimalass} (\Delta\mu)^2 \geq\frac{C \sigma_+^{2}\log d}{n \td_0^2(1 - \td
_0)^2}, \qquad C > 1\quad
\mbox{and} \quad\dd< p/2 -2,
\end{eqnarray}
where $\sigma_+^{2} = \sup_h^* \sigma_h^2$. Then $\widehat{\mathcal
{S}}_d(\alpha_n)$, and hence also $\widehat{\mathcal{S}}_d^c(\alpha
_n)$ are consistent if we either use $\widehat{\sigma}_h^{\min},
\widehat{\sigma}_h^*, \widehat{\sigma}_h^{\diamond},\widehat
{\sigma}_h^{\mathrm{mean}}$, $\widehat{\sigma}_h^-$, $\widehat
{\sigma}_h^+$ or $\widehat{\sigma}_h^{\max}$.
\end{prop}

Note that Proposition~\ref{propconsistentestimationofSd} is only
valid if we use estimators $\widehat{\tau}_h(\td)$ with $\td< \td
_0$. Combining Corollary~\ref{cortime} with Proposition~\ref
{propconsistentestimationofSd} we obtain the following corollary,
which marks the final result of this section.

%co3.7 #&#
\begin{cor}\label{corestimatetimewithinS}
Grant the conditions of Proposition~\ref{propconsistentestimationofSd}. Then
\begin{eqnarray*}
\max_{h \in\widehat{\mathcal{S}}_d^c(\alpha_n)}\bigl\llvert \widehat {\tau}_h(\td) -
\tau_h\bigr\rrvert = \mbox{\scriptsize$\mathcal{O}$}_P(1).
\end{eqnarray*}
\end{cor}

%s4 #&#
\section{Bootstrap under $\mathcal{H}_{\mathcal{A}}$}\label{secpermalt}

As was mentioned before, another option to construct confidence regions
if Theorem~\ref{teoextreme} fails to hold is bootstrapping. In the
context of dependent data, blockwise bootstrap procedures are proposed
in the literature. One of the main problems we face here in this
context is the possibility of change points, and thus a  ``naive'' block
bootstrap can go severely wrong. In the univariate or multivariate
case, possible way outs can be found in \cite{antochhuskova2001} and~\cite{kirch2006}. Here, we employ a different approach that uses the
same idea as in Section~\ref{sectime}, which resulted in consistent
long run variance estimators under both $\mathcal{H}_0$ and $\mathcal
{H}_{A}$. This is outlined in detail in the next section below. In
Section~\ref{secbootII}, we show how this approach can be modified
and simplified. In particular, as a somewhat surprising result, we end
up having a simple  ``naive'' block bootstrap that is consistent.

%s4.1 #&#
\subsection{Bootstrap \textup{I}}\label{secbootI}

Introduce the following notation. For a probability measure $P$ and a
$\sigma$-algebra $\mathcal{G}$, we denote with $P_{| \mathcal{G}}$
the conditional probability with respect to $\mathcal{G}$. Moreover,
we denote with $\mathcal{X} = \sigma ({\mathbf X}_1,\ldots,{\mathbf
X}_n )$ the $\sigma$-algebra generated by the underlying sample.
Pick $0 < \td\leq1/2$, recall that $U_{j,h} = X_{j,h} - \E
[X_{j,h} ]$ and in analogy to $\overline{B}_{n,h}^{\widehat
{\sigma}}$ denote with
\begin{eqnarray*}
\overline{B}_{n,h}^{\widehat{\sigma}^*}(\td) &=& \frac{1}{\widehat
{\sigma}_h^* \sqrt{n}}\sup
_{\td\leq t \leq1 - \td}\biggl\llvert \overline{S}_{n,h}(t) -
\frac{\lceil nt \rceil}{n}\overline {S}_{n,h}(1)\biggr\rrvert, \qquad1 \leq h
\leq d,
\\
\overline{S}_{n,h}(t) &=& \sum_{j = 1}^{\lceil n t \rceil}U_{j,h},
\qquad\overline{T}_{d}^{\widehat{\sigma}^*}(\td) = \max_{1 \leq h
\leq d}
\overline{B}_{n,h}^{\widehat{\sigma}^*}(\td).
\end{eqnarray*}
The objective of this section is to obtain an approximation in the
sense of
\begin{eqnarray*}
\sup_{x \in\R}\bigl\llvert P_{| \mathcal{X}} \bigl(\widehat
{T}_{d,L}^{\widehat{\ss}}(\td) \leq x \bigr) - P \bigl(\overline
{T}_d^{\widehat{\sigma}^*}(\td)\leq x \bigr) \bigr\rrvert = \OO
_P \bigl(n^{-C} \bigr), \qquad C > 0,
\end{eqnarray*}
where\vspace*{1pt} $\widehat{T}_{d,L}^{\widehat{\ss}}(\td)$ is an appropriately
bootstrapped version. To this end, let $K,L$ such that $n = K L$. In
the sequel, $K$ will denote the size of the blocks, and correspondingly
$L$ the number of blocks. For simplicity, we always assume that $K,L
\in\N$, which has no impact on the results. Consider the following
block bounds:
\begin{eqnarray}
\nonumber
\widehat{L}_h^- &=& \sup \bigl\{l \in\N: l K + K/2\leq
\widehat {\tau}_h(\td) n \bigr\},
\nonumber\\[-8pt]\\[-8pt]\nonumber
\widehat{L}_h^+ &=& \inf \bigl\{l \in\N: l K - K/2 \geq\widehat{
\tau}_h(\td) n \bigr\},
\end{eqnarray}
where $\widehat{\tau}_h(\td)$ is as in \eqref{defntauh}. These
estimated limits will allow us to  ``filter'' the contaminated blocks,
and thus allow for a consistent bootstrap procedure. For the actual
construction,
%Note that since $0 < \td\leq1/2$, we always have that $
%\widehat{L}_h^- \thicksim\widehat{L}_h^+ \thicksim L$.
consider the mean estimates
\begin{eqnarray*}
\overline{X}_h^- = \frac{1}{K \widehat{L}_h^-} \sum
_{j = 1}^{K
\widehat{L}_h^-} X_{j,h}, \qquad
\overline{X}_h^+ = \frac{1}{K (L -
\widehat{L}_h^+)} \sum
_{j = K \widehat{L}_h^+ +1}^{K L} X_{j,h},
\end{eqnarray*}
and introduce the random variables
%
%e4.1 #&#
\begin{eqnarray}
\widehat{X}_{j,h} &=& \cases{ X_{j,h} - \overline{X}_h^-,
&\quad if $j \leq K \widehat{L}_h^-$,
\vspace*{3pt}\cr
0, &\quad if $K
\widehat{L}_h^- < j \leq K \widehat{L}_h^+$,
\vspace*{3pt}\cr
X_{j,h} - \overline{X}_h^+, &\quad if $j > K
\widehat{L}_h^+$,}
\end{eqnarray}
and the block variables
%
%e4.2 #&#
\begin{eqnarray}
\widehat{V}_{l,h}(k) = \sum_{j = (l-1)K+1}^{lK}
\widehat {X}_{j,h}\ind(j \leq k), \qquad1 \leq l \leq L, 1 \leq h \leq
d.
\end{eqnarray}
Note the presence of the indicator function $\ind(j \leq k)$ in
$\widehat{V}_{l,h}(k)$, which will allow us to take the maximum within
the individual blocks; see below.

%\begin{comment}
%\begin{eqnarray*}
%\widehat{\mathbf V}_l(k) = \bigl(\widehat{V}_{l,1}(k),\ldots,
%\widehat{V}_{l,d}(k)\bigr)^{\top}
%\end{eqnarray*}
%\end{comment}
Based on $\widehat{V}_{l,h}(\cdot)$, we now have several options for
the construction of a bootstrap, which are among others:
\begin{longlist}[(iii)]
\item[(i)] Multiplier bootstrap.
\item[(ii)] Sampling with replacement.
\item[(iii)] Sampling with no replacement.
\item[(iv)] Mixed versions: (i)${}+{}$(ii) or (i)${}+{}$(iii).
\end{longlist}
In the sequel, we establish results for (i) and (iv), see also~\cite{suppA} for extensions. For $1 \leq l \leq L$, consider the random
variables $\pi(l)$ which take values in the set $\mathcal{L} =
\{1, \ldots, L \}$, and denote with $\pi= \sigma (\pi
(1),\ldots,\pi(L) )$ the corresponding \mbox{$\sigma$-}algebra. The
random variables $\pi(l)$ select the blocks $\widehat{V}_{l,h}(\cdot
)$. Depending on the desired choice of bootstrap, we have that:
\begin{longlist}[(SNR)]
\item[\Mboot]\label{Mboot} $\pi(l) = l$ for $l \in\mathcal{L}$
(deterministic, multiplier bootstrap),
\item[\SRboot]\label{SRboot} $\pi(l)$ are IID and uniformly
distributed over $\mathcal{L}$  (sampling with replacement),
\item[\SNRboot]\label{SNRboot} $\pi(1),\ldots,\pi(L)$ results
from a permutation of $\mathcal{L}$  (sampling with no replacement).
\end{longlist}
Let $\xi_1,\ldots, \xi_L$ be a sequence of IID standard Gaussian
random variables. We then consider the overall statistic $\widehat
{T}_{d,L}^{\widehat{\ss}}(\td)$, defined as
%
%e4.3 #&#
\begin{equation}
\label{defnhatTdL} \widehat{T}_{d,L}^{\widehat{\ss}}(\td) = \max
_{1 \leq h \leq d} \frac{\max_{\lfloor n \td\rfloor\leq k \leq n - \lfloor n \td
\rfloor}}{\widehat{\ss}_{h  \mathcal{X},\pi}\sqrt{n}} \Biggl\llvert \sum
_{l = 1}^L \xi_l \widehat{V}_{\pi(l),h}(k)
- \frac{k}{n}\sum_{l = 1}^L
\xi_l \widehat{V}_{\pi(l),h}(n) \Biggr\rrvert,\hspace*{-30pt}
\end{equation}
where
%
%e4.4 #&#
\begin{eqnarray}
\widehat{\ss}_{h \mid \mathcal{X},\pi}^2 = \frac{1}{K L} \sum
_{l =
1}^L \xi_l^2
\widehat{V}_{\pi(l)}^2(n), \qquad 1 \leq h \leq d,
\end{eqnarray}
denotes the conditional long run variance estimator, which acts as a
replacement for $(\widehat{\sigma}_h^*)^2$. Note that one may also
set $\xi_l^2 =1$ in the definition of $\widehat{\ss}_{h \mid  \mathcal
{X},\pi}^2$. Also note in particular that the maximum (in time) is
taken over $\lfloor n \td\rfloor\leq k \leq n - \lfloor n \td\rfloor
$ in $\widehat{T}_{d,L}^{\widehat{\ss}}(\td)$. Subject to a
specific sample, our bootstrap procedure is now the following.

%al4.1 #&#
\begin{alg}[(Bootstrap algorithm I)]\label{algbootstrap}
\begin{longlist}
\item[\textit{Step} 1.] Pick $0 < \td\leq1/2$, preferably small, compute
$\widehat{\tau}_h(\td)$ for $1 \leq h \leq d$ and select either
\hyperref[Mboot]{\Mboot}, \hyperref[SRboot]{\SRboot} or \hyperref
[SNRboot]{\SNRboot}. Set $m = 1$.

\item[\textit{Step} 2.] Generate $ \{\pi(l) \}_{1 \leq l \leq L}$
according to step~1.\vspace*{1pt}

\item[\textit{Step} 3.] Generate IID $\xi_1,\ldots,\xi_L$ with standard
Gaussian distribution.\vspace*{1pt}

\item[\textit{Step} 4.] Calculate the value of $\widehat{T}_{d,L}^{\widehat
{\ss}}(\td)$ and set $T_m = \widehat{T}_{d,L}^{\widehat{\ss}}(\td)$.

\item[\textit{Step} 5.] Go to step~2 and set $m = m + 1$.
\end{longlist}
\end{alg}

%
%re4.2 #&#
\begin{rem}
As was noted by a reviewer, the definition of $\widehat
{T}_{d,L}^{\widehat{\ss}}(\td)$ in \eqref{defnhatTdL} implies
that both \hyperref[Mboot]{\Mboot} and \hyperref[SNRboot]{\SNRboot}
give an identical procedure. %In~\cite{suppA}, we also discuss
%alternative versions of $\widehat{T}_{d,L}^{\widehat{\ss}}(\td)$ where
%this is not the case with analogue results as given below.
\end{rem}

Stopping Algorithm~\ref{algbootstrap} at $m = M$, we have obtained a
Monte Carlo vector ${\mathbf T}_M =  (T_1,\ldots,T_M )^{\top}$.
For stating consistency results of quantile estimates based on ${\mathbf
T}_M$, it is convenient to parameterize the number of blocks $L$ as $L
= L_n \thicksim n^{\ld}$, where $0 < \ld< 1$. Note that this implies
$K \thicksim n^{1 - \ld}$. The required connection between $\ad,\dd,\ld$ and $p$ is stated in our main assumption below, which can be
considered as mirror conditions for the spatial Assumptions~\ref
{assrelations}.

%as4.3 #&#
\begin{ass}[(Bootstrap assumptions)]\label{assrelationsboot}
Assume that $\ad, \dd, \ld$ and $p > 8$ satisfy:
\begin{longlist}[(B2)]
\item[\Bone]\label{B1} $\dd< \min \{p(1 - \ld)(2\ad-1)(\ad
-1)/\ad- 8\ld, 2\ld(p - 4)  \}/8$,

\item[\Btwo]\label{B2} if $\mathcal{S}_d^c \neq \emptyset$, then for some absolute constant $C > 0$
\begin{eqnarray*}
\bigl\llvert \mathcal{S}_d^c\bigr\rrvert (K \log n
)^{-p/2 + 2} \lesssim n^{-C} \quad\mbox{and}\quad\limsup
_{n \to\infty} \frac{\log n}{K
(\Delta\mu)^2} = 0.
\end{eqnarray*}
\end{longlist}
\end{ass}

We are now ready to state the main result of this section, which
enables us to establish consistency of the above bootstrap procedure.

%th4.4 #&#
\begin{teo}\label{teopureboot}
Grant Assumptions~\ref{assmain},~\ref{assrelationsboot} and one of
\hyperref[Mboot]{\Mboot}, \hyperref[SRboot]{\SRboot} or \hyperref
[SNRboot]{\SNRboot}. Assume in addition that \eqref
{assboundtimechange} holds for some $0 < \td_0 \leq1/2$ if $\mathcal{S}_d^c \neq \emptyset$, and that
\hyperref[S1]{\Sone} is valid. Then for any $0 < \td< \td_0$
\begin{eqnarray*}
\sup_{x \in\R}\bigl\llvert P_{\mid \mathcal{X}} \bigl(\widehat
{T}_{d,L}^{\widehat{\ss}}(\td) \leq x \bigr) - P \bigl(\overline
{T}_d^{\widehat{\sigma}^*}(\td)\leq x \bigr) \bigr\rrvert = \OO
_P \bigl(n^{-C} \bigr), \qquad C > 0.
\end{eqnarray*}
\end{teo}

Let us briefly elaborate on the underlying conditions. Assumption~\ref
{assmain} is our usual temporal and nondegeneracy condition.
\hyperref[B1]{\Bone} and \hyperref[S1]{\Sone} provide the necessary
relation between $\ad,\bd,\dd,\ld$ and the moments $p$. Finally,
\hyperref[B2]{\Btwo} and \eqref{assboundtimechange} are necessary
to control possible change points. Overall, these are rather mild
assumptions. A special point is condition $\td> 0$. It is a purely
technical condition that is required for the proof. One may argue
heuristically that in fact one can also set $\td= 0$ in Theorem~\ref
{teopureboot}, a rigorous argument appears to be rather technical and
lengthy though, and was therefore not pursued. In the simulations in
Section~\ref{secautoempresults}, we set $\td= 0$, and this does
not seem to cause any trouble.

Denote with $\widehat{z}_{\alpha,L}(\td)$ the (conditional)
$1-\alpha$-quantile of $\widehat{T}_{d,L}^{\widehat{\ss}}(\td)$, that is,
\begin{equation}
\widehat{z}_{\alpha,L}(\td) = \inf\bigl\{x\dvtx   P_{|\mathcal{X}}
\bigl(\widehat{T}_{d,L}^{\widehat{\ss}}(\td)\leq x \bigr) \geq 1 - \alpha \bigr\},
\qquad \alpha \in (0,1).
\end{equation}
It then follows from Theorem \ref{teopureboot} that
\[
\bigl|P_{|\mathcal{X}}\bigl(\overline{T}_d^{\diamond,\widehat{\sigma}^*}(\td)\leq \widehat{z}_{\alpha,L}(\td) \bigr)
- \bigl(1 - \alpha\bigr)\bigr| = \OO_p\bigl(n^{-C}), \qquad C > 0,
\]
where $\overline{T}_d^{\diamond,\widehat{\sigma}^*}(\td)$ is a copy of
$\overline{T}_d^{\widehat{\sigma}^*}(\td)$, independent of $\mathcal{X}$.
Standard empirical process theory (cf. \cite{shorak1986}) now implies that we are able to
consistently estimate $\widehat{z}_{\alpha,L}(\td)$ based on $\mathbf{T}_M$ for large
enough $M = M_n$. In particular, in analogy to $\widehat{\mathcal{S}}_{d}^Z(\alpha)$,
the bootstrapped confidence set $\widehat{\mathcal{S}}_{d}(\alpha,L,\td,M)$
can be constructed via
\begin{equation}\label{eqconfregionperm}
\widehat{\mathcal{S}}_{d}(\alpha,L,\td,M) = \bigl\{1 \leq h \leq d\dvtx
 B_{n,h}^{\widehat{\sigma}^*}(\td) \leq \widehat{z}_{\alpha,L}(\td,M) \bigr\},
\end{equation}
where $ \widehat{z}_{\alpha,L}(\td,M)$ is a consistent estimate of
$\widehat{z}_{\alpha,L}(\td)$ based on $\mathbf{T}_M$. Corresponding empirical
examples are given in Section \ref{secsimfactor}.

%s4.2 #&#
\subsection{Bootstrap \textup{II}}\label{secbootII}

Let us step back for a moment to reconsider our original testing,
respectively, estimation problem, that is, construct an estimator for
the stable set $\mathcal{S}_d$. We now make the following observation.
Recall that from the discussion in Section~\ref{sectime} we can
expect that with high probability
\begin{eqnarray*}
\min_{h \in\mathcal{S}_d^c} B_{n,h}^{\widehat{\sigma}^*} > C\sqrt {\log d}
\qquad\mbox{as $n \to\infty$,}
\end{eqnarray*}
for $C > 0$ large enough, if the change in mean is sufficiently strong.
Hence, in order to control the error of estimation for $\widehat
{\mathcal{S}}_d(\alpha)$ or $\widehat{\mathcal{S}}_{d}(\alpha,L,\td)$, we only need to control $\max_{h \in\mathcal{S}_d}
B_{n,h}^{\widehat{\sigma}^*}$. This has interesting consequences for
a bootstrap method, as we will now explain. Recall from Section~\ref
{sectime} that we needed to modify $\widehat{\sigma}_h$ to $\widehat
{\sigma}_h^*$ to avoid the problem of inconsistent variance
estimation. Here, we will actually exploit this problem to our
advantage. More precisely, we construct (conditional) variance
estimators $\widetilde{\ss}_{h|\mathcal{X}}$ that explode for $h \in
\mathcal{S}_d^c$ sufficiently fast, that is, we have with high probability
%
%e4.8 #&#
\begin{eqnarray}
\widetilde{\ss}_{| \mathcal{X}} \gtrsim\sqrt{K}\Delta\mu_h \tau
_h (1 - \tau_h) \qquad\mbox{as $n \to\infty$,}
\end{eqnarray}
where $\Delta\mu_h$ is given in \eqref{defnDeltamu}. Consequently,
we can expect that
\begin{eqnarray*}
\widehat{T}_{d,L}^{\widetilde{\ss}}(\td) = \max_{h \in\mathcal
{S}_d}
\frac{\max_{\lfloor n \td\rfloor\leq k \leq n - \lfloor n
\td\rfloor}}{\widetilde{\ss}_{h \mid \mathcal{X}}\sqrt{n}} \Biggl\llvert \sum_{l = 1}^L
\xi_l \widehat{V}_{l,h}(k) - \frac{k}{n}\sum
_{l =
1}^L \xi_l
\widehat{V}_{l,h}(n) \Biggr\rrvert + \mbox{\scriptsize $
\mathcal{O}$}_P (1 )
\end{eqnarray*}
[see \eqref{eqdefnThattilde} below for $\widehat
{T}_{d,L}^{\widetilde{\ss}}(\td)$]. From Theorem~\ref
{teopureboot}, we then essentially get that
%
%e4.9 #&#
\begin{eqnarray}
\widehat{T}_{d,L}^{\widetilde{\ss}}(\td) \stackrel{d} {=} \max
_{h
\in\mathcal{S}_d} B_{n,h}^{\widehat{\sigma}^*} + \mbox{\scriptsize $
\mathcal{O}$}_P (1 ).
\end{eqnarray}
In other words, $\widehat{T}_{d,L}^{\widetilde{\ss}}(\td)$
automatically \textit{adapts} to the number of unaffected coordinates
and, therefore, allows for a better control of the Type I and II
errors. Note, however, that this will only have a significant impact if
%
%e4.10 #&#
\begin{eqnarray}
\llvert \mathcal{S}_d\rrvert /\bigl\llvert \mathcal{S}_d^c
\bigr\rrvert \to0 \qquad\mbox{as $d \to\infty$.}
\end{eqnarray}
Implementation of this idea will lead to the bootstrap procedure
Algorithm~\ref{algbootstrapII}. However, even more is possible.
Exploiting the explosions another time, we will see that one may
entirely skip estimation of $\widehat{\tau}_h(\td)$, by using a
 ``naive'' bootstrap method. This will lead to Algorithm~\ref
{algbootstrapIII}.

To simplify the following exposition, we only concentrate on multiplier
bootstrap procedures in the remainder of this section. Consider the
overall mean estimator $\overline{X}_h$ and the  ``centered'' random
variables $\widetilde{X}_{j,k}$
%
%e4.11 #&#
\begin{eqnarray}
\overline{X}_h = \frac{1}{n}\sum
_{j = 1}^n X_{j,h}, \qquad
\widetilde{X}_{j,h} = X_{j,h} - \overline{X}_h,
\end{eqnarray}
and the block variables
\begin{eqnarray*}
\widetilde{V}_{l,h}(k) = \sum_{j = (l-1)K+1}^{lK}
\widetilde {X}_{j,h}\ind(j \leq k), \qquad1 \leq l \leq L, 1 \leq h
\leq d.
\end{eqnarray*}
We then construct the (conditional) long run variance estimator
$\widetilde{\ss}_{h\mid \mathcal{X}}^2$ as
%
%e4.12 #&#
\begin{eqnarray}
\widetilde{\ss}_{h\mid \mathcal{X}}^2 = \frac{1}{K L} \sum
_{l = 1}^L \xi_l^2
\widetilde{V}_{l,h}^2(n),
\end{eqnarray}
where $\xi_1,\ldots, \xi_L$ is a sequence of IID standard Gaussian
random variables. Next, pick any $0 < \td\leq1/2$. In analogy to
$\widehat{T}_{d,L}^{\widehat{\ss}}(\td)$, we then consider the
overall statistic $\widehat{T}_{d,L}^{\widetilde{\ss}}(\td)$,
defined as
%
%e4.13 #&#
\begin{eqnarray}
\label{eqdefnThattilde}
\qquad\widehat{T}_{d,L}^{\widetilde{\ss}}(\td) &=& \max
_{1 \leq h \leq d} \frac{\max_{\lfloor n \td\rfloor\leq k \leq n - \lfloor n \td
\rfloor}}{\widetilde{\ss}_{h \mid \mathcal{X}}\sqrt{n}} \Biggl\llvert \sum
_{l = 1}^L \xi_l \widehat{V}_{l,h}(k)
- \frac{k}{n}\sum_{l = 1}^L
\xi_l \widehat{V}_{l,h}(n) \Biggr\rrvert.
\end{eqnarray}
Note that in comparison to $\widehat{T}_{d,L}^{\widehat{\ss}}(\td
)$, we have replaced $\widehat{\ss}_{h \mid \mathcal{X},\pi}$ with
$\widetilde{\ss}_{h \mid \mathcal{X}}$. Subject to a specific sample,
our bootstrap procedure is now the following.

%al4.5 #&#
\begin{alg}[(Bootstrap algorithm II)]\label{algbootstrapII}
\begin{longlist}
\item[\textit{Step} 1.] Pick $0 < \td\leq1/2$, preferably small, compute
$\widehat{\tau}_h(\td)$ for $1 \leq h \leq d$ and $\widetilde{\ss
}_{h\mid  \mathcal{X}}$. Set $m = 1$.

\item[\textit{Step} 2.] Generate IID $\xi_1,\ldots,\xi_L$ with standard
Gaussian distribution.\vspace*{1pt}

\item[\textit{Step} 3.] Calculate the value of $\widehat{T}_{d,L}^{\widetilde
{\ss}}(\td)$ and set $T_m = \widehat{T}_{d,L}^{\widetilde{\ss
}}(\td)$.

\item[\textit{Step} 4.] Go to step~2 and set $m = m + 1$.
\end{longlist}
\end{alg}

By the discussion after Theorem~\ref{teopureboot}, the following
result allows us to establish the consistency of the bootstrap procedure
in Algorithm~\ref{algbootstrapII}.

%th4.6 #&#
\begin{teo}\label{teopurebootII}
Grant Assumptions~\ref{assmain},~\ref{assrelationsboot}. Assume
that \eqref{assboundtimechange} holds for some $0 < \td_0 \leq
1/2$ if $\mathcal{S}_d^c \neq \emptyset$, and that \hyperref[S1]{\Sone} is valid. Then for any $0 < \td<
\td_0$
\begin{eqnarray*}
\sup_{x \in\R}\Bigl\llvert P_{\mid \mathcal{X}} \bigl(\widehat
{T}_{d,L}^{\widetilde{\ss}}(\td) \leq x \bigr) - P \Bigl(\max
_{h
\in\mathcal{S}_d}B_{n,h}^{\widehat{\sigma}^*}(\td)\leq x \Bigr) \Bigr
\rrvert = \mbox{\scriptsize$\mathcal{O}$}_P (1 ).
\end{eqnarray*}
\end{teo}

%$\widetilde{V}_{l,h} = \sum_{j = (l-1)K}^{lK} X_{j,h}$ $
%\overline{X}_{l,h} = 1/K\sum_{j = (l-1)K}^{lK} X_{j,h}$

As a next step, we now show how one can entirely remove the estimation
of $\widehat{\tau}_h(\td)$. The central idea is that one can show
\begin{eqnarray*}
\max_{h \in\mathcal{S}_d^c} \frac{1}{\widetilde{\ss}_{h \mid \mathcal
{X}}\sqrt{n}} \max_{\lfloor n \td\rfloor\leq k \leq n - \lfloor n
\td\rfloor}
\Biggl\llvert \sum_{l = 1}^L
\xi_l \widetilde{V}_{l,h}(k) - \frac{k}{n}\sum
_{l = 1}^L \xi_l
\widetilde{V}_{l,h}(n) \Biggr\rrvert = \OO_P (1 ).
\end{eqnarray*}
Thus, roughly speaking, the  ``explosions'' in both $\widetilde{\ss}_{h
\mid \mathcal{X}}$ and $\xi_l\widetilde{V}_{l,h}(\cdot)$ cancel. Hence,
we automatically obtain the desired relation
\begin{eqnarray*}
\widetilde{T}_{d,L}^{\widetilde{\ss}}(\td) &=& \max_{h \in\mathcal
{S}_d}
\frac{1}{\widetilde{\ss}_{h \mid \mathcal{X}}\sqrt{n}} \max_{\lfloor n \td\rfloor\leq k \leq n - \lfloor n \td\rfloor} \Biggl\llvert \sum
_{l = 1}^L \xi_l \widetilde{V}_{l,h}(k)
- \frac{k}{n}\sum_{l
= 1}^L
\xi_l \widetilde{V}_{l,h}(n) \Biggr\rrvert +
\OO_P (1 )
\\
&\stackrel{d} {=}& \max_{h \in\mathcal{S}_d} B_{n,h}^{\widehat
{\sigma}^*}
+ \OO_P (1 ),
\end{eqnarray*}
where
\begin{eqnarray*}
\widetilde{T}_{d,L}^{\widetilde{\ss}}(\td) &=& \max_{1 \leq h \leq
d}
\frac{1}{\widetilde{\ss}_{h \mid \mathcal{X}}\sqrt{n}} \max_{\lfloor n \td\rfloor\leq k \leq n - \lfloor n \td\rfloor}\Biggl\llvert \sum
_{l = 1}^L \xi_l \widetilde{V}_{l,h}(k)
- \frac{k}{n}\sum_{l
= 1}^L
\xi_l \widetilde{V}_{l,h}(n) \Biggr\rrvert.
\end{eqnarray*}
Note that, in comparison to $\widehat{T}_{d,L}^{\widetilde{\ss}}(\td
)$, we have replaced $\widehat{V}_{l,h}(\cdot)$ with $\widetilde
{V}_{l,h}(\cdot)$. Subject to a specific sample, our bootstrap
procedure is now the following.

%al4.7 #&#
\begin{alg}[(Bootstrap algorithm III)]\label{algbootstrapIII}
\begin{longlist}
\item[\textit{Step} 1.] Pick $0 < \td\leq1/2$, preferably small, compute
$\widetilde{\ss}_{h\mid \mathcal{X}}$ and set $m = 1$.

\item[\textit{Step} 2.] Generate IID $\xi_1,\ldots,\xi_L$ with standard
Gaussian distribution.\vspace*{2pt}

\item[\textit{Step} 3.] Calculate the value of $\widetilde
{T}_{d,L}^{\widetilde{\ss}}(\td)$ and set $T_m = \widetilde
{T}_{d,L}^{\widetilde{\ss}}(\td)$.

\item[\textit{Step} 4.] Go to step~2 and set $m = m + 1$.
\end{longlist}
\end{alg}

As before in Theorems~\ref{teopureboot} and~\ref{teopurebootII},
the following result allows us to conclude consistency of the above
bootstrap procedure.

%th4.8 #&#
\begin{teo}\label{teobootIII}
Grant Assumptions~\ref{assmain},~\ref{assrelationsboot}. Assume
that \eqref{assboundtimechange} holds for some $0 < \td_0 \leq
1/2$ if $\mathcal{S}_d^c \neq \emptyset$, and that \hyperref[S1]{\Sone} is valid. If in addition we have
with
\begin{equation}\label{eqthmbootIIInondegcondi}
\liminf_{d \to \infty} P\Bigl(\max_{h \in
 \mathcal{S}_d}B_{n,h}^{{\sigma}}(\td_0) \geq x_d \Bigr) = 1
\end{equation}
for some monotone increasing sequence $x_d \to \infty$,
then for any $0 < \td < \td_0$
\begin{eqnarray*}
\sup_{x \in\R}\Bigl\llvert P_{| \mathcal{X}} \bigl(\widetilde
{T}_{d,L}^{\widetilde{\ss}}(\td) \leq x \bigr) - P \Bigl(\max
_{h
\in\mathcal{S}_d}B_{n,h}^{\widehat{\sigma}^*}(\td)\leq x \Bigr) \Bigr
\rrvert = \mbox{\scriptsize$\mathcal{O}$}_P (1 ).
\end{eqnarray*}
\end{teo}

%
%re4.9 #&#
\begin{rem}\label{remthmbootIIInondeg}
Assumption \eqref{eqthmbootIIInondegcondi} is a mild
nondegeneracy condition, and is only violated in the extreme case
where $\lim_{d \to\infty}\max_{h \in\mathcal
{S}_d}B_{n,h}^{{\sigma}}(\td_0) = \OO_P (1 )$.
\end{rem}

%s4.3 #&#
\subsection{Discussion of bootstrap procedures}\label{secdiscussbootstrap}

In the previous sections, we have seen that all three bootstrap
Algorithms~\ref{algbootstrap},~\ref{algbootstrapII} and~\ref
{algbootstrapIII} are consistent alternatives to Theorem~\ref
{teoextreme} and Proposition~\ref{propsimpleboot}. In particular,
they do not require any assumptions on the spatial dependence
structure. On the other hand, there are also some deficits that we will
briefly outline.
\begin{description}
\item[Computational cost:] Particularly if $d$ gets larger, the
computational costs and time become a relevant issue.
\item[Homogeneity:] All bootstrap procedures require global blocks in
order to reflect the underlying dependence structure. This in turn
requires a certain homogeneity of temporal dependence of the data.
\item[Sensitivity:] As simulations reveal, the number and thus size of
the blocks may have a huge impact, and in some cases the results appear
to be rather sensitive in this respect, and there is also an interplay
with the required homogeneity, mentioned above. This problem of
blocklength selection is already well known in the literature in the
univariate or multivariate case; see, for instance,~\cite
{lahiribook2003}. A~simple problematic example is given in
Section~\ref{secsimfactor}, Table~\ref{tabbootgoeswrong}.
\item[Large \textit{d} small \textit{n}:] If $d$ is rather large compared to $n$,
one should take $L$ as large as possible to avoid or at least weaken
some of the above problems. In particular, one should keep in mind that
one multiplies and thus  ``models'' the time series with only $L$ i.i.d.
Gaussian random variables. However, a large $L$ results in a small $K$,
and thus a possible failure in capturing the temporal dependence via
the blocks.
\end{description}

From these considerations, a bootstrap procedure is only recommended if
the dimension $d$ is not too large compared to the sample size $n$ ($d
\leq n$ appears to still yield good results), and if the vast majority
of the data can be expected to be homogenous (a few outliers do not
hurt). Otherwise, the parametric bootstrap depicted in Proposition~\ref
{propsimpleboot} is recommended.

%s5 #&#
\section{Examples}\label{secexamples}

In this section, we discuss some prominent and leading examples from
the literature that fit into our framework. To keep the exposition at a
reasonable length, our main focus lies on ARMA$(\mathfrak{p},\mathfrak
{q})$ and GARCH$(\mathfrak{p},\mathfrak{q})$ models, but our setup
also contains many more nonlinear time series, as will be briefly
discussed. We mainly focus on examples that fulfill Assumptions~\ref
{assmain} and~\ref{assrelations}. Of course, this implies that these
are also valid examples for the bootstrap procedures. An important
example are factor models in Example~\ref{exfactor}, which highlight the
usefulness of bootstrap procedures.

In the setting of Theorem~\ref{teoextreme}, the spatial decay
condition \hyperref[S3]{\Sthree} plays a key role.
The (multivariate) time series literature contains a huge variety of
process that meet \hyperref[A1]{\Tone}. Especially for nonlinear
time series such as GARCH-models, iterated random functions and the
like, we refer to~\cite{fanbook2005,hormann2009,poetscher1999}
and the many references there. We thus concentrate on giving examples
for $\{{\mathbf X}_k\}_{k \in\Z}$, where \hyperref[S3]{\Sthree} holds.
More precisely, we give examples for two parameter processes $ \{
X_{k,h} \}_{k,h \in\Z}$ where the key conditions \hyperref
[A1]{\Tone} and \hyperref[S3]{\Sthree} are valid.

We recall the following convention. Throughout this section, $0 < C <
\infty$ denotes an arbitrary, absolute constant that may vary from
line to line.

%
%ex5.1 #&#
\begin{ex}[(Linear processes)]\label{exlinhigh}
A common way to model multivariate linear models with finite dimension
$d$ is by
%
%e5.1 #&#
\begin{eqnarray}
\label{eqexlin} {\mathbf X}_k = \sum_{l = 0}^{\infty}
{\mathbf R}_l {\mathbf Z}_{k-l},
\end{eqnarray}
where $ \{{\mathbf R}_i \}_{i \in\N}$ is a sequence of $d\times
d$ matrices, and $ \{{\mathbf Z}_{k} \}_{k \in\Z}$ is a
sequence of $d$-dimensional vectors, usually i.i.d. However,
describing (weak) spatial dependence in this model when $d$ is large is
not at all straightforward, even if one assumes a simple spatial
structure for ${\mathbf Z}_{k}$, for example, a linear process. In
addition, using high-dimensional matrices for modelling purposes is
only advisable if the matrices are sparse. The problem of transferring
multivariate linear models, in particular the autoregressive
multivariate setup to a high-dimensional setting is currently a very
active field of research, particularly in connection with panel data or
factor models. For example, in~\cite{chudik2011}, various sparsity
constraints are discussed to introduce the IVAR (infinite-dimensional
vector-autoregression). Other approaches are offered in~\cite
{chentwosample2010,ChoFryzlewicz2012Preprint}. Here, we will first
follow the approach taken in the latter, before coming back to \eqref
{eqexlin}. Let $ \{\varepsilon_{k,h} \}_{k,h \in\Z}$ be a
sequence such that $\varepsilon_k =  \{\varepsilon_{k,h} \}_{h
\in\Z}$ is i.i.d. for $k \in\Z$. We then introduce the high-dimensional MA($\infty$, $\varepsilon$) process as
\begin{eqnarray*}
X_{k,h} = \sum_{i = 0}^{\infty}
\alpha_{i,h} \varepsilon_{k-i,h}\qquad\mbox{for $k,h \in\Z$ and $
\alpha_{i,h} \in\R$.}
\end{eqnarray*}
Naturally, we require some conditions on the numbers $\alpha_{i,h}$ to
guarantee its existence. We do this in one sweep, by stating conditions
such that Assumptions \hyperref[A1]{\Tone} and \hyperref[S3]{\Sthree
} are valid in addition.

%pr5.2 #&#
\begin{prop}\label{proplinearY}
Suppose that $\llvert  \gamma_{i,j}^{\varepsilon}\rrvert   = \llvert  \E
[\varepsilon_{k,i}\varepsilon_{k,j} ]\rrvert   \leq C \log
(\llvert  i-j\rrvert  )^{-2-\delta}$ for $\llvert  i-j\rrvert   \geq2$ and $\delta> 0$. If \hyperref
[A2]{\Ttwo} holds and also
\begin{eqnarray*}
\sup_h^*\llvert \alpha_{i,h}\rrvert \lesssim
i^{-\ud}\qquad\mbox {with }\ud> 5/2,
\end{eqnarray*}
then $ \{X_{k,h} \}_{k,h \in\Z}$ meets Assumptions
\textup{\hyperref[A1]{\Tone}} and \textup{\hyperref[S3]{\Sthree}}.
\end{prop}

As a special case, we may now consider ARMA($\mathfrak{p},\mathfrak
{q}$, $\varepsilon$) processes, which we introduce as
\begin{eqnarray*}
X_{k,h} = \alpha_{1,h}^* \varepsilon_{k,h} + \cdots+
\alpha_{\mathfrak
{p},h}^* \varepsilon_{k - \mathfrak{p},h} + \beta_{1,h}^*
X_{k-1,h} + \cdots+ \beta_{\mathfrak{q},h}^* X_{k - \mathfrak{q},h},
\end{eqnarray*}
$\alpha_{1,h}^*,\ldots,\alpha_{\mathfrak{p},h}^*, \beta_{1,h}^*,\ldots,\beta_{\mathfrak{q},h}^* \in\R$. As in the univariate case, we
consider the polynomials
%
%e5.2 #&#
\begin{eqnarray}
\label{eqrelpolys} {\mathbf A}_h(z) = \sum_{j = 0}^{\mathfrak{p}}
\alpha_{j,h}^* z^{j}, \qquad{\mathbf B}_h(z) = \sum
_{j = 0}^{\mathfrak{q}} \beta_{j,h}^*
z^{j},
\end{eqnarray}
where ${\mathbf A}_h(z)$ and ${\mathbf B}_h(z)$ are assumed to be relative
prime. Then following, for instance, \cite{timeseriesbrockwell}, one
readily verifies the following result.

%pr5.3 #&#
\begin{prop}
If the associated polynomials ${\mathbf C}_h(z) = {\mathbf A}_h(z) {\mathbf
B}_h^{-1}(z)$ satisfy $\inf_h^*\llvert  {\mathbf C}_h(z)\rrvert   > 0$ for
$\llvert  z\rrvert   \leq1$, then $X_{k,h}$ admits a causal representation
\begin{eqnarray*}
X_{k,h} = \sum_{i = 0}^{\infty}
\alpha_{i,h} \varepsilon_{k-i,h}\qquad \mbox{where }\sup
_h^*\llvert \alpha_{k,h}\rrvert \lesssim
q^{k}\mbox{ for }0 < q < 1.
\end{eqnarray*}
\end{prop}

It is now easy to see that we may choose $\ad> 5/2$ arbitrarily large,
hence Assumption \hyperref[A1]{\Tone} holds. The validity of
\hyperref[S3]{\Sthree} can be obtained as in Proposition~\ref
{proplinearY}. Next, we demonstrate how model \eqref{eqexlin} fits
into our framework. Recall that ${\mathbf Z}_{k} =  \{Z_{k,h} \}
_{h \in\Z}$. We impose the following conditions.

%as5.4 #&#
\begin{ass}\label{asslinear}
The sequence $ \{{\mathbf Z}_{k} \}_{k \in\Z}$ is IID, and for
$p > 4$:
\begin{longlist}[(ii)]
\item[(i)] $\gamma_{i,j}^{Z} = \E [Z_{k,i} Z_{k,j}  ] \leq
C\log(\llvert  i-j\rrvert   + 2)^{-2 - \delta}$, $\delta> 0$,
\item[(ii)] ${\mathbf R}_l =  (r_{i,j}^{(l)} )_{1 \leq i,j \leq
d}$ with $\llvert  r_{i,j}^{(l)}\rrvert   \leq C (l+1)^{-\qd} (\llvert  i-j\rrvert   +
1)^{-\pd}$, $\qd,\pd> 2$.
\end{longlist}
\end{ass}

Condition (ii) is mild, allowing for a large variety of matrix
sequences ${\mathbf R}_l$. We have the following result.

%pr5.5 #&#
\begin{prop}\label{proplingen}
Assume that Assumptions~\ref{asslinear} and \textup{\hyperref[A2]{\Ttwo}}
hold. Then $ \{X_{k,h} \}_{k,h \in\Z}$ meets Assumptions
\textup{\hyperref[A1]{\Tone}} and \hyperref[S3]{\Sthree}.
\end{prop}

Based on the above proposition, one can derive a related result for
multivariate ARMA processes, we omit the details.
\end{ex}

%
%ex5.6 #&#
\begin{ex}[(Factor models)]\label{exfactor}
In econometric theory, it is often believed that the dynamics of a
multivariate or high-dimensional time series ${\mathbf X}_k$ can be
described via so-called (normally unobserved) common factors ${\mathbf Z}_k
\in\R^{d'}$, where it is usually assumed in the literature that $d'$
is much smaller than $d$. This amounts to the model
%
%e5.3 #&#
\begin{eqnarray}
\label{eqmodelfact} {\mathbf X}_k = {\mathbf R} {\mathbf Z}_k +
\xiv_k, \qquad k \in\Z,
\end{eqnarray}
where ${\mathbf R} =  (r_{i,j} )_{1 \leq i \leq d, 1 \leq j \leq
d'}$ is a $d \times d'$ matrix of factor loadings, and $\xiv_k = \{\xi
_{k,h}\}_{h \in\Z}$ denotes the noise sequence. We also denote with
$\ss_{h,\xi}^2$ the coordinate-wise standard deviation of $\{\xi
_{k,h}\}_{k \in\Z}$, and with $\phi_{k,i,j}^{Z} = \E
[Z_{k,i}Z_{0,j} ]$, $\phi_{k,i,j}^{\xi} = \E [\xi
_{k,i}\xi_{0,j} ]$. We then make the following assumptions.

%as5.7 #&#
\begin{ass}\label{assfactor}
For $\delta> 0$ and $p > 4$, we have:
\begin{longlist}[(iii)]
\item[(i)] $ \{{\mathbf Z}_{k} \}_{k \in\Z}$ and $ \{\xiv
_{k} \}_{k \in\Z}$ are independent and both satisfy Assumption
\hyperref[A1]{\Tone},

\item[(ii)] $\phi_{k,i,j}^{Z},\phi_{k,i,j}^{\xi}\leq C (k+1)^{-\qd
} (\llvert  i-j\rrvert   + 1)^{-\pd}$, $\qd,\pd> 2$ and $\inf_h^* \ss_{h,\xi}^2 > 0$,

\item[(iii)] $\sup_{i}^*\sum_{j = 1}^{d'} \llvert  r_{i,j}\rrvert   < \infty$ and
$\llvert  \sum_{j = 1}^{d'} r_{i_1,j}r_{i_2,j}\rrvert   \leq C  (\log
\llvert  i_1 - i_2\rrvert   )^{-2 -\delta}$ for
$\llvert  i_1 - i_2\rrvert   \geq2$.
\end{longlist}
\end{ass}

The above assumptions are related to those of Assumption~\ref
{asslinear}. This comes as no surprise, since both process are very
similar. As we shall see, rather straightforward computations show that
the corresponding spatial correlation matrix $\Sigma_d =  (\rho
_{i,j} )_{1 \leq i,j \leq d}$ only needs to satisfy
%
%e5.4 #&#
\begin{eqnarray}
\label{eqconditrhoopt} \llvert \rho_{i,j}\rrvert \leq C \bigl(\log\llvert i - j
\rrvert \bigr)^{-2 -\delta}\quad \delta> 0\qquad\mbox{if }\llvert i - j\rrvert
\geq2.
\end{eqnarray}
We now have the following result.

%pr5.8 #&#
\begin{prop}\label{propfactor}
Assume that $d' = d'_n$ with $d'_n \to\infty$ and $d' \leq d$.
Suppose that Assumptions~\ref{assfactor}, \hyperref[S1]{\Sone} and
\hyperref[S2]{\Stwo} hold. Then Theorem \ref{teoextreme} is valid.
\end{prop}

The above result shows that under reasonable assumptions, high-dimen\-sional
factor models fit into our framework. Note that the limit
distribution is pivotal, no additional information on ${\mathbf R}$ is
required. This is important from a statistical point of view, since the
factor loadings ${\mathbf R}$ are usually unobservable in practice. In this
context, the question arises whether a structural condition like
%
%e5.5 #&#
\begin{eqnarray}
\label{eqfactordecaycondi} \Biggl\llvert \sum_{j = 1}^{d'}
r_{i_1,j}r_{i_2,j}\Biggr\rrvert \lesssim \bigl(\log\llvert
i_1 - i_2\rrvert \bigr)^{-2 -\delta}\qquad\mbox{if }
\llvert i_1 - i_2\rrvert \geq2
\end{eqnarray}
is necessary to obtain a pivotal limit distribution. If \eqref
{eqfactordecaycondi} does not hold, one can still show via Theorems \ref{teoextreme}, \ref{teoupperbound} and the triangle
inequality that with probability one
%
%e5.6 #&#
\begin{eqnarray}
\label{equpperlower} \qquad\liminf_{n} \frac{\max_{1 \leq h \leq d} B_{n,h}^{\widehat{\sigma
}}}{\sqrt{\log d}} > 0 \quad
\mbox{and}\quad\limsup_{n} \frac{\max_{1 \leq h \leq d} B_{n,h}^{\widehat{\sigma}}}{\sqrt{\log d}} \leq
\frac{1}{\sqrt{2}},
\end{eqnarray}
hence $\sqrt{\log d}$ is the right scaling, even without \eqref
{eqfactordecaycondi}. However, determining the exact limit
distribution in the absence of \eqref{eqfactordecaycondi} seems to
be very difficult, and is likely to depend on ${\mathbf R}$, questioning
its usefulness for applications. In fact, if we drop condition \eqref
{eqfactordecaycondi} in Assumption~\ref{assfactor}(iii), a
pivotal result like Theorem \ref{teoextreme} cannot hold as the
next result shows.

%pr5.9 #&#
\begin{prop}\label{propfactorno}
Assume that the conditions of Proposition~\ref{propfactor} hold, with
the exception that we do not have \eqref{eqfactordecaycondi}. Then
universal sequences $a_d$, $b_d$, only depending on $d$ such that
\begin{eqnarray*}
a_d \Bigl(\max_{1 \leq h \leq d} B_{n,h}^{\widehat{\sigma}}
- b_d \Bigr)
\end{eqnarray*}
converges in distribution to a nondegenerate limit do not exist.
\end{prop}

Proposition~\ref{propfactorno} tells us that an exact fluctuation
control without any intrinsic knowledge on ${\mathbf R}$ is not possible.
In this sense, Assumption~\ref{assfactor} seems to be near the
minimum requirements to obtain a pivotal, nonparametric result like
Theorem~\ref{teoextreme}. In any case, relation \eqref
{equpperlower} tells us that we always remain in control of the Type
I error, and the possible loss in power is only marginal.
\end{ex}

%
%ex5.10 #&#
\begin{ex}[(GARCH process)]\label{exgarchhigh}
In this example, we discuss one possible way to extend the constant
conditional GARCH model (CCG) of Bollerslev~\cite{bollerslev}. If the
dimension $d$ is fixed, related multivariate extensions can be found in
the literature; see, for instance,~\cite{aueetal2009}. Here, we
define the GARCH$(\mathfrak{p},\mathfrak{q}, \varepsilon)$ sequence as
\begin{eqnarray*}
X_{k,h} &=& \varepsilon_{k,h} \ss_{k,h}\qquad
\mbox{where } \{\ss _{k,h} \}_{k,h \in\Z}\mbox{ meets},
\\
\ss_{k,h}^2 &=& \eta_h + \alpha_{1,h}
\ss_{k - 1,h}^2 + \cdots+ \alpha_{\mathfrak{p},h}
\ss_{k - \mathfrak{p},h}^2 + \beta_{1,h} X_{k - 1,h}^2
+ \cdots+ \beta_{\mathfrak{q},h} X_{k - \mathfrak{q},h}^2,
\end{eqnarray*}
with $\eta_h, \alpha_{1,h},\ldots,\alpha_{\mathfrak{p},h}, \beta
_{1,h},\ldots,\beta_{\mathfrak{q},h} \in\R$. Note that $\mathfrak
{p}$ and $\mathfrak{q}$ denote the maximal degree of $\alpha_{i,h}$,
$\beta_{i,h}$. Possible undefined $\alpha_{i,h}$ and $\beta_{i,h}$
are replaced with zeros. As in the univariate case, a crucial quantity
in this context is
%
%e5.7 #&#
\begin{eqnarray}
\gamma_C = \max_{1 \leq h \leq d}\sum
_{i = 1}^{r} \bigl\llVert \alpha _{i,h} +
\beta_{i,h} \varepsilon_{i,h}^2\bigr\rrVert
_{p/2}\qquad\mbox{with } r = \max\{\mathfrak{p},\mathfrak{q}\}.
\end{eqnarray}
If $\gamma_C < 1$, then it can be shown that $ \{X_{k,h} \}
_{k,h \in\Z}$ is stationary (cf.~\cite{bougerol1992aop}). We have
the following result, establishing a link between the underlying
parameters and Assumption~\ref{assmain}.

%pr5.11 #&#
\begin{prop}\label{propextendgarch}
Suppose that $\gamma_{C} < 1$ and $\gamma_{i,j}^{\varepsilon} = \E
 [\varepsilon_{k,i} \varepsilon_{k,j}  ]$ satisfies
\begin{eqnarray*}
\bigl\llvert \gamma_{i,j}^{\varepsilon}\bigr\rrvert \leq C \bigl(\log
\llvert i - j\rrvert \bigr)^{-2 -\delta} \delta> 0\qquad\mbox{if }\llvert i -
j\rrvert \geq2.
\end{eqnarray*}
Then $ \{X_{k,h} \}_{k,\in\Z, h \in\N}$ meets Assumptions
\hyperref[A1]{\Tone} and \hyperref[S3]{\Sthree}.
\end{prop}
\end{ex}

In~\cite{suppA}, we additionally discuss time series that arise as
iterated random functions. Moreover, as in the univariate case, many
more examples can be constructed based on the vast time series
literature (cf.~\cite{fanbook2005,hormann2009,poetscher1999}).
Also note that any combination of the previous examples fulfills
Assumption~\ref{assmain}. This means that in one coordinate we may
have a GARCH model, but in another coordinate, the process has a linear
dynamic, and so on.

%s6 #&#
\section{Empirical results and applications}\label{secempresults}

In the empirical part of the paper, we first discuss the implications
and relevance of our assumptions for real data sets. We then move on to
the computation of critical values. In the third part, we asses the
accuracy and behavior of $\widehat{\mathcal{S}}_{d}^c$ in a small
simulation study. In~\cite{suppA}, the S\&P 500
companies over a period of one year,
with a particular emphasis of detecting companies with an unusual behavior.
%s6.1 #&#
\subsection{Assumptions and real data}\label{secassplusdata}

The necessary assumptions of Theorem~\ref{teoextreme} can be divided
into temporal and spatial conditions. Assumption~\ref{assmain}
concerns temporal dependence, and is standard in the literature (cf.~\cite{wu2005}). We therefore focus on the spatial conditions, in
particular \hyperref[S3]{\Sthree}. This condition implicitly assumes
that the coordinates of the data-vector ${\mathbf X}_k$ can be ordered such
that two coordinates $X_{k,i}$ and $X_{k,j}$ become less dependent as
the difference $\llvert  i-j\rrvert  $ increases. In many cases, the data at hand
already has a natural ordering with corresponding weak spatial
dependence. Such examples can be found in the ever growing literature
on high-dimensional covariance estimation (cf.~\cite{cai2011}),
where spatial dependence is modelled (or expressed) by a banding or
block-wise structure of the matrix. Note that in this case, the order
of the coordinates is essential for the covariance estimator and needs
to be specified in advance. In our setup, however, the advantageous
order need not be known explicitly to the practitioner due to the
maximum statistic.

%\begin{comment}
%In fact, in light of Remark~\ref{remassm}, we only require that the
%data is (or can be) ordered in such a way that it consists of clusters
%$\CC_h$, $h = 1,2,\ldots$, where the clustersize $\left\vert \CC_h\right\vert $
%itself may grow in $n$. Within each cluster, the data may be highly
%dependent, as long as $\left\vert \rho_{i,j}\right\vert  < 1$ for $i \neq j$ where $i,j \in
%\bigl\{1,\ldots,d\bigr\} \setminus{\mathcal T}$ and ${\mathcal T}
%\subset\bigl\{1,\ldots,d\bigr\}$ with cardinality $\oo\bigl(d\bigr)$
%(cf. Remark~\ref{remremoveindexset}). The spatial condition
%\hyperref[A1]{\Tone} now effectively concerns the clusters $\CC_h$,
%and requires that $\CC_i, \CC_j$ become less dependent as $\left\vert i-j\right\vert $ gets
%large, possibly very slowly, see Example~\ref{exfactor}. Note that
%this ``clustering'' is actually also the basic idea of the proof, and
%is also connected to properties of the Gaussian distribution.
%\end{comment}

If a spatial condition like \hyperref[S3]{\Sthree} cannot be assumed
to hold (see Example~\ref{exfactor}), we can use the bootstrap
procedures from Section~\ref{secpermalt}. However, at least some
preliminary considerations should be made; see the short discussion in
Section~\ref{secdiscussbootstrap}. One way to check whether the
permutation bootstrap is necessary is by means of a PCA. The literature
on factor models provides a simple heuristic test (cf.~\cite
{connor1993}) in this direction. Compute the largest empirical
eigenvalue $\widehat{\lambda}_1$ of the empirical correlation matrix
$\widehat{\Sigma}_d$. In the presence of a common factor, $\widehat
{\lambda}_1$ will explode with rate $d$, that is,
%
%e6.1 #&#
\begin{eqnarray}
\liminf_{d \to\infty} \widehat{\lambda}_1/d \geq C > 0
\qquad\mbox{a.s.}
\end{eqnarray}
Hence, if $\widehat{\lambda}_1/d$ is small, a common factor is rather
unlikely or its overall influence very weak, and a bootstrap is not
necessary. As a final remark, let us mention that if controlling the
Type I error is essential, the parametric bootstrap is highly
recommended as a first tool for inference. The empirical results
regarding the bootstrap in Section~\ref{secsimfactor} reveal that
the behavior may be significantly influenced by the choice of the
number of blocks $L$ and the connected size $K$, which makes
controlling the Type I error not so easy.\looseness=1

%s6.2 #&#
\subsection{Critical values}\label{seccriticalvalues}

%
%t1 #&#
\begin{table}[b]
\tabcolsep=0pt
\caption{Parametric bootstrap. Sample size $n \in\{100,250,500\}$,
dimension $d \in\{100,250,500\}$}\label{tabcritvalues}
\begin{tabular*}{\tablewidth}{@{\extracolsep{\fill}}@{}lcccccccccccc@{}}
\hline
& \multicolumn{3}{c}{$\bolds{n = 100}$} & \multicolumn{3}{c}{$\bolds{n = 250}$}
& \multicolumn{3}{c}{$\bolds{n = 500}$} & \multicolumn{3}{c@{}}{\textbf{Numerical}}
\\[-6pt]
& \multicolumn{3}{c}{\hrulefill} & \multicolumn{3}{c}{\hrulefill}
& \multicolumn{3}{c}{\hrulefill} & \multicolumn{3}{c@{}}{\hrulefill}
\\
$\bolds{d}$ & $\bolds{100}$ & $\bolds{250}$ & $\bolds{500}$ & $\bolds{100}$ & $\bolds{250}$ & $\bolds{500}$ & $\bolds{100}$ & $\bolds{250}$ &
$\bolds{500}$ & $\bolds{100}$ & $\bolds{250}$ & $\bolds{500}$\\
\hline
${\mathbf q_{0.9}}$ & 1.83 & 1.93 & 2.00 & 1.88 & 1.99 & 2.07 & 1.9 & 2.02 & 2.10 & 1.95 & 2.05 & 2.14\\
${\mathbf q_{0.95}}$ & 1.91 & 2.00 & 2.10 & 1.97 & 2.07 & 2.15 & 1.99 & 2.10 & 2.19 & 2.03 & 2.15 & 2.22\\
${\mathbf q_{0.975}}$ & 1.98 & 2.07 & 2.15 & 2.05 & 2.15 & 2.22 & 2.08 & 2.19 & 2.28 & 2.12 & 2.21 & 2.30\\
${\mathbf q_{0.99 }}$ & 2.07 & 2.17 & 2.24 & 2.15 & 2.25 & 2.31 & 2.19 & 2.30 & 2.36 & 2.22 & 2.32 & 2.40\\
\hline
\end{tabular*}
\end{table}

Deriving reasonable critical values for extreme value statistics is a
delicate issue. The root of the problem typically lies in the slow
convergence rate of extreme value statistics. In our case, the domain
of attraction is the Gumbel distribution, and the rate of convergence
(for Gaussian random variables) is no better than $\OO (\log
n^{-1} )$; see~\cite{hall1979}. Hence, using the normalizing
sequences $e_d, f_d$ given in Theorem~\ref{teoextreme} may not be the
best thing to do. On the other hand, as is demonstrated by Proposition
\ref{propsimpleboot}, approximative critical values can either be
obtained by a parametric bootstrap or numerical computations. In
principle, there are two methods for obtaining critical values in case
of the parametric bootstrap.
\begin{longlist}[(a)]
\item[(a)] Simulate $\max_{1 \leq h \leq d} B_{n,h}^Z$ directly.
\item[(b)] Estimate $F_n(x) = P (B_{n,h}^Z \leq x )$, and
obtain the critical values via $1 - \alpha= F_n(z_{\alpha})^{d}$.
\end{longlist}
Method (b) is more flexible and was used to obtain the results. The
corresponding critical values are tabulated in Table~\ref
{tabcritvalues}. A total of $10^6$ MC-runs was used to compute each
critical value. Generally speaking, the quantiles obtained by numerical
computations (Table~\ref{tabcritvalues}, column  ``Numerical'') are
larger. This can be explained by the fact that in the ``discrete''
version $B_{n,h}^{Z}$, the maximum is taken over the set $\{1,\ldots,n\}$, whereas in the limit $\B_h$, the supremum is taken over the
whole interval $[0,1]$, which is a larger set, and hence leads to this
relation. In case of the multiplier bootstrap, very similar results
are obtained in the same Gaussian setup. Empirical evidence for the
validity of the multiplier bootstrap in the presence of dependence and
change points is given in Section~\ref{secsimfactor}, where critical
values are tabulated in Table~\ref{tabcritvaluesfactor}.
%\begin{comment}
%The code for generating all quantiles (for both bootstrap methods) is
%available at~\cite{r-code}, and all tables can exactly be reproduced
%using the corresponding seed.
%\end{comment}

%s6.3 #&#
\subsection{Simulation study}

In this subsection, we investigate the Type I error and power of the
estimator $\widehat{\mathcal{S}}_d^c$ in a small simulation study. We
consider estimates originating from the parametric as well as the
permutation bootstrap. To assess the power, several alternatives are
considered: we insert artificial changes in certain coordinates $h$ at
$\tau_h \in \{(2i+1)/10 \}$, $0 \leq i \leq4$ with size
$\delta/10$ where $\delta\in \{0,0.25,0.5,0.75,1 \}$. We
then study the behavior and estimation accuracy on the sets
\begin{eqnarray*}
\mathcal{S}_d^c = \mathcal{S}_{d,1}^c
\uplus\mathcal{S}_{d,2}^c \uplus\mathcal{S}_{d,3}^c
\uplus\mathcal{S}_{d,4}^c \uplus\mathcal
{S}_{d,5}^c,
\end{eqnarray*}
where
\begin{eqnarray*}
\mathcal{S}_{d,i}^c = \bigl\{h \in\mathcal{S}_d^c:
\tau_h \in \bigl[(i-1)/5, i/5\bigr) \bigr\}, \qquad1 \leq i \leq5.
\end{eqnarray*}
Note that this means that we check whether the test detects a change
and also classifies the time of change correctly. As a measure of
comparison, we evaluate the relative estimation accuracy (in $\%$) as
\begin{eqnarray*}
r_{d,i}^c = \frac{\E [\llvert  \widehat{\mathcal{S}}_{d,i}^c
\cap\mathcal{S}_{d,i}^c \rrvert   ]}{\llvert  \mathcal
{S}_{d,i}^c\rrvert  } \times100, \qquad1 \leq i
\leq5,
\end{eqnarray*}
where the mean $\E [\llvert  \widehat{\mathcal{S}}_{d,i}^c \cap
\mathcal{S}_d^c \rrvert   ]$ is estimated from the overall
simulated sample. This gives an accurate measure of the performance of
the test procedure. We also consider the coordinatewise Type I error,
described by the probability
\begin{eqnarray*}
\mathrm{TI}_{h} = P \bigl(h \in\widehat{\mathcal{S}}_{d}^c
\cap \mathcal{S}_d \bigr), \qquad h \in\mathcal{S}_d.
\end{eqnarray*}
Note that if $ \{X_{k,h} \}_{k \in\Z, h \in\N}$ is a
stationary random field (which is the case in all of our simulations),
then $\mathrm{TI}_{h}$ is the same for all $h \in\mathcal{S}_d$ and
can be written as
\begin{eqnarray*}
\mathrm{TI}_{h} = \frac{\E [\llvert  \widehat{\mathcal
{S}}_{d}^c \setminus\mathcal{S}_{d}^c \rrvert   ]}{\llvert  \mathcal
{S}_d\rrvert  }, \qquad h \in
\mathcal{S}_d.
\end{eqnarray*}
To allow for reproducibility and transparency, all simulations use
exactly the same random seed, and also the sets $\mathcal{S}_{d,i}^c$
remain the same. This implies that for $n,d$ fixed, the Type I error
$\mathrm{TI}_{h}$ remains \textit{the same} for all $\delta$.
Natural exceptions are only when $\delta= 0$ or the long run
simulations in Tables~\ref{tabbootII} and~\ref{tabbootIII}
concerning the bootstrap results. The number of change points for each
$\mathcal{S}_{d,i}^c$ is set to $10$ for $d = 100$ and $15$ for $d =
250$. This gives a total amount of changes $\llvert  \mathcal
{S}_{100}^c\rrvert   = 50$ and $\llvert  \mathcal{S}_{250}^c\rrvert   = 75$.
As sample size, we consider $n \in\{100,250\}$ and 1000 MC runs for
each setting, unless stated otherwise. We use two different models for
our simulations, namely autoregressive and factor models. In case of
the factor model, we also investigate the behavior of the bootstrap
Algorithms~\ref{algbootstrapII} and~\ref{algbootstrapIII}.

%
%t2 #&#
\begin{table}[t]
\tabcolsep=0pt
\caption{Sample size $n = 100$, dimension $d = 100$, $\mathrm
{TI}_{h}^* = \mathrm{TI}_{h} \times100$, $\alpha= 0.05$, $\widehat
{\sigma}_h^*$}\label{tabbootII}
\begin{tabular*}{\tablewidth}{@{\extracolsep{\fill}}@{}ld{1.1}d{2.1}d{2.1}d{2.1}d{1.1}d{1.2}d{1.1}d{2.1}d{2.1}d{2.1}d{1.1}d{1.2}@{}}
\hline
& \multicolumn{6}{c}{\textbf{Bootstrap II} $\bolds{d = 100}$} & \multicolumn{6}{c@{}}{\textbf{Parametric} $\bolds{d = 100}$}\\[-6pt]
& \multicolumn{6}{c}{\hrulefill} & \multicolumn{6}{c@{}}{\hrulefill}
\\
$\bolds{\delta}$ & \multicolumn{1}{c}{$\bolds{r_{d,1}^c}$}
& \multicolumn{1}{c}{$\bolds{r_{d,2}^c}$}
& \multicolumn{1}{c}{$\bolds{r_{d,3}^c}$}
& \multicolumn{1}{c}{$\bolds{r_{d,4}^c}$}
& \multicolumn{1}{c}{$\bolds{r_{d,5}^c}$}
& \multicolumn{1}{c}{$\bolds{\mathrm{TI}_{h}^*}$}
& \multicolumn{1}{c}{$\bolds{r_{d,1}^c}$}
& \multicolumn{1}{c}{$\bolds{r_{d,2}^c}$}
& \multicolumn{1}{c}{$\bolds{r_{d,3}^c}$}
& \multicolumn{1}{c}{$\bolds{r_{d,4}^c}$}
& \multicolumn{1}{c}{$\bolds{r_{d,5}^c}$}
& \multicolumn{1}{c@{}}{$\bolds{\mathrm{TI}_{h}^*}$}\\
\hline
$ 0$ & \multicolumn{1}{c}{--} & \multicolumn{1}{c}{--} & \multicolumn{1}{c}{--} & \multicolumn{1}{c}{--} & \multicolumn{1}{c}{--} & 3.12 & \multicolumn{1}{c}{--} & \multicolumn{1}{c}{--} & \multicolumn{1}{c}{--} & \multicolumn{1}{c}{--} & \multicolumn{1}{c}{--} & 1.86\\
$0.025$ & 0.3 & 3.3 & 8.7 & 3.2 & 0.2 & 3.2 & 0.3 & 2.6 & 6.6 & 1.8 & 0.1 & 1.86\\
$0.05$ & 1.3 & 18.6 & 34.6 & 18.1 & 1.0 & 3.52 & 0.8 & 12.8 & 26.4 & 11.6 & 0.6 & 1.86\\
$0.075$ & 3.6 & 48.2 & 68.4 & 46.5 & 3.3 & 3.76 & 2.7 & 39.3 & 59.3 & 37.2 & 1.6 & 1.86\\
$0.1$ & 8.7 & 72.0 & 85.9 & 70.1 & 6.2 & 2.9 & 6.2 & 65.8 & 82.6 & 64.3 & 5.0 & 1.86\\
\hline
\end{tabular*}
\end{table}

%
%t3 #&#
\begin{table}[b]
\tabcolsep=0pt
\caption{Sample size $n = 100$, dimension $d = 100$, $\mathrm
{TI}_{h}^* = \mathrm{TI}_{h} \times100$, $\alpha= 0.05$, $\widehat
{\sigma}_h^*$}\label{tabbootIII}
\begin{tabular*}{\tablewidth}{@{\extracolsep{\fill}}@{}ld{1.1}d{2.1}d{2.1}d{2.1}d{1.1}d{1.2}d{1.1}d{2.1}d{2.1}d{2.1}d{1.1}d{1.2}@{}}
\hline
& \multicolumn{6}{c}{\textbf{Bootstrap III} $\bolds{d = 100}$} & \multicolumn{6}{c@{}}{\textbf{Parametric} $\bolds{d = 100}$}\\[-6pt]
& \multicolumn{6}{c}{\hrulefill} & \multicolumn{6}{c@{}}{\hrulefill}
\\
$\bolds{\delta}$ & \multicolumn{1}{c}{$\bolds{r_{d,1}^c}$}
& \multicolumn{1}{c}{$\bolds{r_{d,2}^c}$}
& \multicolumn{1}{c}{$\bolds{r_{d,3}^c}$}
& \multicolumn{1}{c}{$\bolds{r_{d,4}^c}$}
& \multicolumn{1}{c}{$\bolds{r_{d,5}^c}$}
& \multicolumn{1}{c}{$\bolds{\mathrm{TI}_{h}^*}$}
& \multicolumn{1}{c}{$\bolds{r_{d,1}^c}$}
& \multicolumn{1}{c}{$\bolds{r_{d,2}^c}$}
& \multicolumn{1}{c}{$\bolds{r_{d,3}^c}$}
& \multicolumn{1}{c}{$\bolds{r_{d,4}^c}$}
& \multicolumn{1}{c}{$\bolds{r_{d,5}^c}$}
& \multicolumn{1}{c@{}}{$\bolds{\mathrm{TI}_{h}^*}$}\\
\hline
$ 0$ & \multicolumn{1}{c}{--} & \multicolumn{1}{c}{--} & \multicolumn{1}{c}{--} & \multicolumn{1}{c}{--} & \multicolumn{1}{c}{--} & 2.3 & \multicolumn{1}{c}{--} & \multicolumn{1}{c}{--} & \multicolumn{1}{c}{--} & \multicolumn{1}{c}{--} & \multicolumn{1}{c}{--} & 1.86\\
$0.025$ & 0.3 & 2.9 & 7.3 & 2.5 & 0.1 & 2.4 & 0.3 & 2.6 & 6.6 & 1.8 & 0.1 & 1.86\\
$0.05$ & 0.8 & 15.2 & 30.0 & 14.0 & 0.8 & 2.52 & 0.8 & 12.8 & 26.4 & 11.6 & 0.6 & 1.86\\
$0.075$ & 3.0 & 43.5 & 63.0 & 40.7 & 2.2 & 2.66 & 2.7 & 39.3 & 59.3 & 37.2 & 1.6 & 1.86\\
$0.1$ & 9.3 & 70.5 & 85.0 & 68.8 & 5.9 & 2.92 & 6.2 & 65.8 & 82.6 & 64.3 & 5.0 & 1.86\\
\hline
\end{tabular*}
\end{table}

%
%t4 #&#
\begin{table}[t]
\tabcolsep=0pt
\caption{Sample size $n = 100$, dimension $d \in\{100,250\}$,
$\mathrm{TI}_{h}^* = \mathrm{TI}_{h} \times100$, $\alpha= 0.05$,
$\widehat{\sigma}_h^*$}\label{tabaltvaluesAR100}
\begin{tabular*}{\tablewidth}{@{\extracolsep{\fill}}@{}ld{1.2}d{2.2}d{2.2}d{2.2}d{1.2}d{1.2}d{1.2}d{1.2}d{2.2}d{2.2}d{1.2}d{1.2}@{}}
\hline
& \multicolumn{6}{c}{\textbf{Parametric} $\bolds{d = 100}$} & \multicolumn{6}{c@{}}{\textbf{Parametric} $\bolds{d = 250}$}\\[-6pt]
& \multicolumn{6}{c}{\hrulefill} & \multicolumn{6}{c@{}}{\hrulefill}
\\
$\bolds{\delta}$ & \multicolumn{1}{c}{$\bolds{r_{d,1}^c}$}
& \multicolumn{1}{c}{$\bolds{r_{d,2}^c}$}
& \multicolumn{1}{c}{$\bolds{r_{d,3}^c}$}
& \multicolumn{1}{c}{$\bolds{r_{d,4}^c}$}
& \multicolumn{1}{c}{$\bolds{r_{d,5}^c}$}
& \multicolumn{1}{c}{$\bolds{\mathrm{TI}_{h}^*}$}
& \multicolumn{1}{c}{$\bolds{r_{d,1}^c}$}
& \multicolumn{1}{c}{$\bolds{r_{d,2}^c}$}
& \multicolumn{1}{c}{$\bolds{r_{d,3}^c}$}
& \multicolumn{1}{c}{$\bolds{r_{d,4}^c}$}
& \multicolumn{1}{c}{$\bolds{r_{d,5}^c}$}
& \multicolumn{1}{c@{}}{$\bolds{\mathrm{TI}_{h}^*}$}\\
\hline
$ 0$ & \multicolumn{1}{c}{--} & \multicolumn{1}{c}{--} & \multicolumn{1}{c}{--} & \multicolumn{1}{c}{--} & \multicolumn{1}{c}{--} & 2.01 & \multicolumn{1}{c}{--} & \multicolumn{1}{c}{--} & \multicolumn{1}{c}{--} & \multicolumn{1}{c}{--} & \multicolumn{1}{c}{--} & 1.23\\
$0.025$ & 0.19 & 2.46 & 6.67 & 2.24 & 0.11 & 2.01 & 0.06 & 1.43 & 5.05 & 1.76 & 0.06 & 1.23\\
$0.05$ & 0.74 & 13.1 & 28.1 & 12.5 & 0.55 & 2.01 & 0.37 & 9.97 & 22.8 & 9.87 & 0.27 & 1.23\\
$0.075$ & 2.54 & 38.8 & 58.4 & 38.0 & 1.85 & 2.01 & 1.52 & 32.3 & 51.6 & 30.7 & 0.95 & 1.23\\
$0.1$ & 7.21 & 67.2 & 82.5 & 65.8 & 4.61 & 2.01 & 5.03 & 61.2 & 77.3 & 59.8 & 2.90 & 1.23\\
\hline
\end{tabular*}
\end{table}

%
%t5 #&#
\begin{table}[b]
\tabcolsep=0pt
\caption{Sample size $n = 250$, dimension $d \in\{100,250\}$,
$\mathrm{TI}_{h}^* = \mathrm{TI}_{h} \times100$, $\alpha= 0.05$,
$\widehat{\sigma}_h^*$}\label{tabaltvaluesAR250}
\begin{tabular*}{\tablewidth}{@{\extracolsep{\fill}}@{}ld{1.2}d{2.2}d{2.1}d{2.2}d{1.2}d{1.2}d{2.2}d{2.2}d{2.1}d{2.2}d{2.2}d{1.2}@{}}
\hline
& \multicolumn{6}{c}{\textbf{Parametric} $\bolds{d = 100}$} & \multicolumn{6}{c@{}}{\textbf{Parametric} $\bolds{d = 250}$}\\[-6pt]
& \multicolumn{6}{c}{\hrulefill} & \multicolumn{6}{c@{}}{\hrulefill}
\\
$\bolds{\delta}$ & \multicolumn{1}{c}{$\bolds{r_{d,1}^c}$}
& \multicolumn{1}{c}{$\bolds{r_{d,2}^c}$}
& \multicolumn{1}{c}{$\bolds{r_{d,3}^c}$}
& \multicolumn{1}{c}{$\bolds{r_{d,4}^c}$}
& \multicolumn{1}{c}{$\bolds{r_{d,5}^c}$}
& \multicolumn{1}{c}{$\bolds{\mathrm{TI}_{h}^*}$}
& \multicolumn{1}{c}{$\bolds{r_{d,1}^c}$}
& \multicolumn{1}{c}{$\bolds{r_{d,2}^c}$}
& \multicolumn{1}{c}{$\bolds{r_{d,3}^c}$}
& \multicolumn{1}{c}{$\bolds{r_{d,4}^c}$}
& \multicolumn{1}{c}{$\bolds{r_{d,5}^c}$}
& \multicolumn{1}{c@{}}{$\bolds{\mathrm{TI}_{h}^*}$}\\
\hline
$ 0$ & \multicolumn{1}{c}{--} & \multicolumn{1}{c}{--} & \multicolumn{1}{c}{--} & \multicolumn{1}{c}{--} & \multicolumn{1}{c}{--} & 0.94 & \multicolumn{1}{c}{--} & \multicolumn{1}{c}{--} & \multicolumn{1}{c}{--} & \multicolumn{1}{c}{--} & \multicolumn{1}{c}{--} & 0.50\\
$0.025$ & 0.04 & 4.96 & 13.9 & 5.02 & 0.04 & 0.94 & 0.02 & 3.22 & 10.2 & 3.52 & 0.03 & 0.50\\
$0.05$ & 0.76 & 42.5 & 65.8 & 41.5 & 0.56 & 0.94 & 0.25 & 34.1 & 58.7 & 33.4 & 0.27 & 0.50\\
$0.075$ & 5.35 & 84.7 & 95.5 & 83.9 & 4.08 & 0.94 & 1.98 & 79.6 & 93.6 & 79.9 & 2.28 & 0.50\\
$0.1$ & 19.7 & 96.3 & 99.6 & 95.6 & 16.9 & 0.94 & 10.7 & 96.2 & 99.6 & 95.6 & 11.4 & 0.50\\
\hline
\end{tabular*}
\end{table}

%s6.3.1 #&#
\subsubsection{Autoregressive models}\label{secautoempresults}

We use the following model. We take $Y_{k,h}$ as an MA(100) process
%
%e6.2 #&#
\begin{eqnarray}
\label{eqsimMA} \qquad Y_{k,h} = \sum_{i = 0}^{99}
\alpha_i \varepsilon_{k,h-i}, \qquad\alpha _i = 0.1
\llvert i\rrvert ^{-3}\mbox{ and }\varepsilon_{k,h} \sim\Gaussian
\bigl(0,\ss^2 \bigr), \ss= 0.1.
\end{eqnarray}
We then consider the ARMA$(2,2)$ model
%
%e6.3 #&#
\begin{eqnarray}
\label{eqsimARMA} \qquad X_{k,h} = 0.2 X_{k-1,h} - 0.3 X_{k-2,h}
-0.1 Y_{k,h} + 0.2 Y_{k-1,h}, \qquad1 \leq h \leq d.
\end{eqnarray}
Note that we stick to the same model in each coordinate, which makes
the comparison and analysis easier and more transparent. Throughout
this section, the nominal level of all tests is $\alpha= 0.05$, that
is, we always use the corresponding quantiles ${\mathbf q_{0.95}}$. We
first analyze the parametric bootstrap. The corresponding results are
given in Tables~\ref{tabaltvaluesAR100} and~\ref
{tabaltvaluesAR250}. Note that in both tables, the row with $\delta
= 0$ corresponds to the empirical levels of the test. The Type I error
is slightly different from the cases where $\delta> 0$ (not visible
due to rounded values), which is due to the fact that $\mathcal{S}_d =
 \{1,\ldots,d \}$ if $\delta= 0$, and $\mathcal{S}_d
\subset \{1,\ldots,d \}$ otherwise. Observe that small
changes are found with difficulty if the sample size is small, and this
effect naturally gets amplified in higher dimensions. The power for
bigger samples/changes is however very reasonable. As expected, the
test loses power if the time of change $\tau_h$ moves away from the
center $1/2$. Unreported simulations exhibit a similar behavior in case
of GARCH-sequences, or tests for changes in the second moment or variance.

Next, we briefly discuss a possible effect in the choice of variance
estimator. In the previous results, estimator $\widehat{\sigma}_h^*$
was used; see Section~\ref{sectime} to recall the definition. As one
comparison, we now use $\widehat{\sigma}_h^{\diamond}$;
see Table~\ref{tabsigdiamondvaluesAR250} for corresponding results. An
interesting phenomena appears. We note that $\widehat{\sigma}_h^*$
yields the better results if $\tau_h = 1/2$, and $\widehat{\sigma
}_h^{\diamond}$ if the change is more away from $1/2$. This is a
little surprising, since one can show that for large enough $n$,
$\widehat{\sigma}_h^*$ has the smaller MSE. A possible explanation
could be the quality of estimation of~$\widehat{\tau}_h$, and the
actual choice of $B_{\tau}$.

We now turn to the behavior of the bootstrap procedures, more
precisely, we consider Algorithms~\ref{algbootstrapII} and~\ref
{algbootstrapIII}, where we ``illegally'' set $\td= 0$. We use the
same model \eqref{eqsimARMA}. In order to obtain an overall feasible
computational time, we restrict ourselves to the setup where $n = 100$,
$d = 100$ and we only used $100$ overall simulations for comparison
(note: comparing the parametric results indicates that there actually
is not much difference between 100 or 1000 simulations). Moreover, we
only use $M = 100$ Monte Carlo runs for the bootstrap procedures.
Arguably, this might be too low to obtain a necessary accuracy for a
$95\%$ quantile, but it turns out that this is not the case. We choose
$L = 25$ as the number of blocks, and thus $K = 4$ for the block
length. The results of Algorithm~\ref{algbootstrapII} are given in
Table~\ref{tabbootII}. Even though we only set $M = 100$, we get
slightly better results than the parametric procedure. Observe that the
results are also conservative. The behavior of Algorithm~\ref
{algbootstrapIII} in Table~\ref{tabbootIII} is slightly worse, but
overall very similar.

%s6.3.2 #&#
\subsubsection{Factor models and number of block length effect}\label{secsimfactor}

%
%t6 #&#
\begin{table}[t]
\tabcolsep=0pt
\caption{Sample size $n = 250$, dimension $d \in\{100,250\}$,
$\mathrm{TI}_{h}^* = \mathrm{TI}_{h} \times100$, $\alpha= 0.05$,
$\widehat{\sigma}_h^{\diamond}$}\label{tabsigdiamondvaluesAR250}
\begin{tabular*}{\tablewidth}{@{\extracolsep{\fill}}@{}ld{2.2}d{2.2}d{2.1}d{2.2}d{2.2}d{1.2}d{2.2}d{2.2}d{2.1}d{2.2}d{2.2}d{1.2}@{}}
\hline
& \multicolumn{6}{c}{\textbf{Parametric} $\bolds{d = 100}$} & \multicolumn{6}{c@{}}{\textbf{Parametric} $\bolds{d = 250}$}\\[-6pt]
& \multicolumn{6}{c}{\hrulefill} & \multicolumn{6}{c@{}}{\hrulefill}
\\
$\bolds{\delta}$ & \multicolumn{1}{c}{$\bolds{r_{d,1}^c}$}
& \multicolumn{1}{c}{$\bolds{r_{d,2}^c}$}
& \multicolumn{1}{c}{$\bolds{r_{d,3}^c}$}
& \multicolumn{1}{c}{$\bolds{r_{d,4}^c}$}
& \multicolumn{1}{c}{$\bolds{r_{d,5}^c}$}
& \multicolumn{1}{c}{$\bolds{\mathrm{TI}_{h}^*}$}
& \multicolumn{1}{c}{$\bolds{r_{d,1}^c}$}
& \multicolumn{1}{c}{$\bolds{r_{d,2}^c}$}
& \multicolumn{1}{c}{$\bolds{r_{d,3}^c}$}
& \multicolumn{1}{c}{$\bolds{r_{d,4}^c}$}
& \multicolumn{1}{c}{$\bolds{r_{d,5}^c}$}
& \multicolumn{1}{c@{}}{$\bolds{\mathrm{TI}_{h}^*}$}\\
\hline
$ 0$ & \multicolumn{1}{c}{--} & \multicolumn{1}{c}{--} & \multicolumn{1}{c}{--} & \multicolumn{1}{c}{--} & \multicolumn{1}{c}{--} & 1.51 & \multicolumn{1}{c}{--} & \multicolumn{1}{c}{--} & \multicolumn{1}{c}{--} & \multicolumn{1}{c}{--} & \multicolumn{1}{c}{--} & 0.99\\
$0.025$ & 0.12 & 8.12 & 16.6 & 8.64 & 0.06 & 1.42 & 0.04 & 6.25 & 12.6 & 6.47 & 0.05 & 1.01\\
$0.05$ & 1.58 & 53.7 & 64.1 & 54.0 & 1.26 & 1.42 & 0.80 & 48.0 & 57.5 & 48.4 & 0.68 & 1.01\\
$0.075$ & 9.12 & 88.6 & 90.2 & 88.6 & 8.56 & 1.42 & 5.95 & 87.5 & 85.6 & 86.8 & 5.93 & 1.01\\
$0.1$ & 29.5 & 96.8 & 97.9 & 96.4 & 29.4 & 1.42 & 23.5 & 96.2 & 95.6 & 96.1 & 23.2 & 1.01\\
\hline
\end{tabular*}
\end{table}

%
%t7 #&#
\begin{table}[b]
\tabcolsep=0pt
\caption{Bootstrap Alg. \textup{II}, \textup{III}. Sample size $n = 250$, dimension $d =
100$, $\delta= 0$, $\alpha_F \in\{0.1,0.3\}$}\label{tabcritvaluesfactor}
\begin{tabular*}{\tablewidth}{@{\extracolsep{\fill}}@{}lcccccccc@{}}
\hline
& \multicolumn{4}{c}{\textbf{Bootstrap} $\bolds{\alpha_F = 0.1}$}
& \multicolumn{4}{c@{}}{\textbf{Bootstrap} $\bolds{\alpha_F = 0.3}$}\\[-6pt]
& \multicolumn{4}{c}{\hrulefill}
& \multicolumn{4}{c@{}}{\hrulefill}
\\
& \multicolumn{2}{c}{$\bolds{n = 250}$} & \multicolumn{2}{c}{$\bolds{n = 250}$} & \multicolumn{2}{c}{$\bolds{n = 250}$} & \multicolumn{2}{c@{}}{$\bolds{n = 250}$}
\\
$\bolds{K \times L}$ & \multicolumn{2}{c}{$\bolds{5\times50}$} & \multicolumn{2}{c}{$\bolds{10 \times25}$} & \multicolumn{2}{c}{$\bolds{5\times50}$} & \multicolumn{2}{c@{}}{$\bolds{10 \times25}$} \\
$\bolds{d}$ & \multicolumn{2}{c}{$\bolds{100}$} & \multicolumn{2}{c}{$\bolds{100}$} & \multicolumn{2}{c}{$\bolds{100}$} & \multicolumn{2}{c@{}}{$\bolds{100}$}\\[-6pt]
\multicolumn{1}{@{}c}{\hrulefill} & \multicolumn{2}{c}{\hrulefill} & \multicolumn{2}{c}{\hrulefill} & \multicolumn{2}{c}{\hrulefill} & \multicolumn{2}{c@{}}{\hrulefill}\\
\textbf{Algorithm} & \textbf{II} & \textbf{III} & \textbf{II} & \textbf{III} & \textbf{II} & \textbf{III} & \textbf{II} & \textbf{III}\\
\hline
${\mathbf q_{0.9}}$ & 1.69 & 1.73 & 1.69 & 1.77 & 1.47 & 1.49 & 1.44 & 1.51\\
${\mathbf q_{0.95}}$ & 1.75 & 1.83 & 1.78 & 1.87 & 1.56 & 1.59 & 1.56 & 1.61\\
${\mathbf q_{0.975}}$ & 1.81 & 1.89 & 1.91 & 1.94 & 1.69 & 1.68 & 1.63 & 1.75\\
${\mathbf q_{0.99 }}$ & 1.89 & 1.95 & 2.05 & 2.14 & 1.82 & 1.80 & 1.73 & 1.85\\
\hline
\end{tabular*}
\end{table}

We consider a factor model that shows that block and parametric
bootstrap may behave very differently. As explained in Example~\ref
{exfactor}, this is the case if overall dependence on certain factors
is present. To allow for a comparison to the autoregressive model, we
use the same model for the general dynamics. We take $Y_{k,h}$ as an
MA(100) process
%
%e6.4 #&#
\begin{eqnarray}
\label{eqsimMA2}
\qquad Y_{k,h} = \sum_{i = 0}^{99}
\alpha_i \varepsilon_{k,h-i}, \qquad\alpha _i = 0.1
\llvert i\rrvert ^{-3}\mbox{ and }\varepsilon_{k,h} \sim\Gaussian
\bigl(0,\ss^2 \bigr), \ss= 0.1.
\end{eqnarray}
We then consider the ARMA$(2,2)$-factor model
\begin{eqnarray*}
X_{k,h} = \alpha_F F_{k} + 0.2
X_{k-1,h} - 0.3 X_{k-2,h} -0.1 Y_{k,h} + 0.2
Y_{k-1,h}, \qquad1 \leq h \leq d,
\end{eqnarray*}
where $\alpha_F \geq0$ is a constant, and the factors $ \{
F_k \}_{k \in\Z}$ are IID standard Gaussian random variables.
The primary focus in this section is to demonstrate the effect of high
spatial dependence and the size $L$ on the quantiles. The pronounced
effect of the factors is visible in Table~\ref
{tabcritvaluesfactor}, where the critical values of Algorithms~\ref
{algbootstrapII} and~\ref{algbootstrapIII} are tabulated for
$\alpha_F \in\{0.1,0.3\}$. A value of $\delta= 0$ and $M = 1000$
were used in the simulations. We see that the factor has an expected
significant reciprocal impact on the quantiles, that is, larger factors
result in lower quantiles. We also observe that in this setup, the
number of blocks $L \in\{25,50\}$ does not have a notable impact.
A~slight outlier seems to be the results of Algorithm~\ref
{algbootstrapII}(II) in case of $L = 25$. Also note that
particularly the results about the more extremal quantiles ${\mathbf
q_{0.975}}$ and ${\mathbf q_{0.99}}$ have to be considered with care,
since ``only'' $M = 1000$ was used.

Unreported simulations show that the power and size are different from
the autoregressive model. Particularly if $\alpha_F$ is large, (e.g.,
$\alpha_F = 0.3$), one has to consider a larger size of change $\delta
> 0.1$ %\in\bigl\{0.1,0.15,0.2,0.25 \bigr\}$
in order to obtain visible effects. The reason for this loss in power
is the (considerably) larger long-run variance $\sigma_h$ in this
model, induced by the factor loading $\alpha_F$. More precisely, since
we scale by a (larger) consistent estimate of $\sigma_h$, changes
become harder to detect [see also \eqref{eqelemtarylowerboundchange}].

Finally, we take a look at Table~\ref{tabbootgoeswrong}, which
reveals that the choice of $L$ may have a serious impact. Here, we set
$\alpha_F = 0$ to allow for a comparison to the results in
Sections~\ref{seccriticalvalues} and~\ref{secautoempresults}. We
observe that raising $K$ only by one from $K = 4$ to $K = 5$ leads to
much larger quantiles. In view of the results presented in Section~\ref
{secautoempresults}, these would lead to a high loss in power. The
setup itself appears rather harmless; we note, however, that $d/n = 1$
have the same size, unlike to the situation in Table~\ref
{tabcritvaluesfactor} where $d/n = 2/5$. It appears that at least if
$d/n \geq1$, block bootstrap procedures based on multipliers can
require careful tuning. Particularly if $d \gg n$, the parametric
bootstrap seems to be the more stable option.

%
%t8 #&#
\begin{table}
\tabcolsep=0pt
\caption{Bootstrap Alg. \textup{II}, \textup{III}. Sample size $n = 100$, dimension $d =
100$, $\delta= 0$, $\alpha_F = 0$}\label{tabbootgoeswrong}
\begin{tabular*}{\tablewidth}{@{\extracolsep{\fill}}@{}lcccc@{}}
\hline
& \multicolumn{2}{c}{$\bolds{n = 100}$} & \multicolumn{2}{c@{}}{$\bolds{n = 100}$} \\
$\bolds{K \times L}$ & \multicolumn{2}{c}{$\bolds{5\times20}$} & \multicolumn{2}{c@{}}{$\bolds{4 \times25}$}\\
$\bolds{d}$ & \multicolumn{2}{c}{$\bolds{100}$} & \multicolumn{2}{c@{}}{$\bolds{100}$}\\[-6pt]
\multicolumn{1}{@{}c}{\hrulefill}  & \multicolumn{2}{c}{\hrulefill} & \multicolumn{2}{c@{}}{\hrulefill}\\
\textbf{Algorithm} & \textbf{II} & \textbf{III} & \textbf{II} & \textbf{III} \\
\hline
${\mathbf q_{0.9}}$ & 1.91 & 2.04 & 1.75 & 1.81 \\
${\mathbf q_{0.95}}$ & 2.08 & 2.23 & 1.89 & 1.92 \\
${\mathbf q_{0.975}}$ & 2.25 & 2.44 & 1.95 & 2.02 \\
${\mathbf q_{0.99 }}$ & 2.38 & 2.68 & 2.05 & 2.11 \\
\hline
\end{tabular*}
\end{table}

%s7 #&#
\section{Proofs}\label{secproofsofextremebridge}

All proofs together with additional results are given in detail in~\cite{suppA}.

%\begin{appendix}
%\section{}
%\end{appendix}

\section*{Acknowledgements}
I am indebted to the Associate Editor and the reviewers for a careful
reading of the manuscript, and the many thoughtful remarks and comments
that considerable helped to improve the results and presentation. I
warmly thank Wei Biao Wu for his comments on the proofs, and Lajos
Horv\'{a}th for his advice. I also thank Holger Drees for his expertise
on extremal quantile estimation. Special thanks to Josua G\"{o}smann for
pointing out some errors.

\begin{supplement}[id=suppA]
\stitle{Supplement to  ``Uniform change point tests in high dimension.''}
\slink[doi]{10.1214/15-AOS1347SUPP} %[doi,text={...}] - jei reikia
%suskaldyti doi
\sdatatype{.pdf}
\sfilename{AOS1347\_supp.pdf}
\sdescription{The supplement contains all proofs.
In addition, a detailed analysis of a data set (S\&P 500) is presented.}
\end{supplement}

% imsref loaded by linak, 2015-07-09 08:55:35

\printaddresses
\end{document}